%% file: main.tex
\documentclass[11pt]{article}

\usepackage[margin=1in]{geometry}
\usepackage[round]{natbib}
\usepackage[utf8]{inputenc}
\usepackage{amssymb,amsthm,mathtools}
\usepackage{booktabs}
\usepackage{graphicx}
\usepackage{xcolor}
\usepackage{microtype}
\usepackage[hidelinks]{hyperref}

\newtheorem{theorem}{Theorem}
\newtheorem{proposition}[theorem]{Proposition}
\newtheorem{lemma}[theorem]{Lemma}
\newtheorem{corollary}[theorem]{Corollary}
\theoremstyle{remark}
\newtheorem{remark}[theorem]{Remark}

\newcommand{\KL}{\mathrm{KL}}
\newcommand{\TV}{\mathrm{TV}}
\newcommand{\N}{\mathcal N}
\newcommand{\E}{\mathbb E}
\newcommand{\R}{\mathbb R}
\newcommand{\Beta}{\mathrm{Beta}}
\newcommand{\nhard}{n_{\mathrm{hard}}}
\newcommand{\ngauss}{n_{\mathrm{Gauss}}}

\begin{document}

\title{When Is Membership Graded? Identifiability and Resolution Limits of a Graded Gaussian Mixture Model}
\author{Hendrata Dharmawan\qquad Peter Chin\\[4pt]
Dartmouth College}
\date{}
\maketitle

\begin{abstract}
Mixed-membership models describe each unit as a graded blend of latent types. Before the blend is interpreted, the data should be able to rule out the simpler story in which every unit belongs to a single type. We ask when they can, in a scalar graded Gaussian mixture model, in which each observation is a Beta-weighted convex combination of independent Gaussian draws from two types, the mechanism behind functional mixed-membership models and the normal compositional model. The model is identifiable, unequal variances included, but the identifying information sits in the far tails of the marginal. For the symmetric model we give a closed-form hard-membership mixture that matches the graded marginal in five cumulants, and a Le Cam bound turns it into a certified sample size below which no test separates graded from hard membership with summed error under one half. The threshold is governed by the concentration $c$ relative to $c_0\approx d^2$, where $d$ is the type half-separation in within-type standard deviations. At $d=3$ it is about 2,600 observations at $c=3$ and 77,000 at $c=6$. Most of the difficulty comes from learning $c$. On 100 simulated datasets with graded membership, a correctly specified posterior under the default priors of a published Bayesian sampler for this model class reports hard membership every time, a weak prior gives intervals for the concentration spanning orders of magnitude, and the sampler's own answers depend on where it starts.
\end{abstract}

\input{sections/intro}
\input{sections/model}
\input{sections/identifiability}
\input{sections/graded_vs_hard}
\input{sections/gaussian_end}
\input{sections/fixed_c}
\input{sections/inference}
\input{sections/discussion}

\section*{Use of generative AI}
In this work, we used generative AI tools for formulating precise mathematical claims consistent with the standard terminology used in the literature, writing the proof and double-checking notations and consistencies in the proof, identifying gaps in the proof, designing and providing feedback for the suite of experiments as well as suggesting additional ablations, helping analyze the results and identifying bugs in the code that produced such results, and identifying typos / non-material algebraic errors in the proof.

We have not used generative AI tools for supplying the main mathematical ideas of the proof (all ideas including key mathematical ideas in the proof belonged to the authors), result interpretations, main coding tasks (other than identifying bugs).

Additionally, tasks such as cleaning and reformatting datasets, assisting with translation are not applicable to this work.

Additionally, we used generative AI tools to create plotting and graphing code suitable for latex given experimental results, suggesting the suites of experiments and the range of parameters to run them over, drafting parts of the paper as well as grammar correction, performing literature search in order to discover other works that have been done, summarizing other papers related to our works, analyzing our paper versus those other papers in order to ensure we do not overclaim novelty, brainstorming next experiment ideas given those other works, grammar and style correction.

We have reviewed all AI-assisted work. We take responsibility for the final content of this work, including text, claims, or artifacts produced with the aid of generative AI.

\bibliographystyle{plainnat}
\bibliography{refs}

\clearpage
\appendix
\input{sections/appendix}

\end{document}

%% file: sections/intro.tex
\section{Introduction}\label{sec:intro}

Mixed-membership models give each unit a vector of weights over latent types in place of a single label. For discrete data the idea is standard \citep{pritchard2000,blei2003,erosheva2004,airoldi2008}. For continuous measurements, a common construction writes each observation as a convex combination of independent draws, one from each type. This is the mechanism of the functional and multivariate mixed-membership models of \citet{marco2024functional,marco2024flexible} and of the normal compositional model used in hyperspectral unmixing \citep{stein2003,eismann2007,eches2010}. In both settings the weights are the quantity of scientific interest. In unmixing they are the abundances of materials within a pixel. In the EEG studies of \citet{marco2024functional,marco2024flexible}, a child's membership measures how far a peak in the alpha band has emerged above an aperiodic $1/f$ background.

Before asking where a unit lies between two types, one can ask whether units lie between types at all. The alternative is that every unit belongs to exactly one type and the population is an ordinary two-component mixture. We ask when data can tell these two stories apart, and find that for many parameter values they cannot at the sample sizes one usually has. We work in the simplest model in which the question is well posed: a scalar observation
\begin{equation}\label{eq:model}
X = U\,Y_1 + (1-U)\,Y_2,\qquad Y_k\sim\N(\mu_k,\sigma_k^2),\quad U\sim\Beta(a,b),
\end{equation}
with $U$, $Y_1$ and $Y_2$ independent. We call this the graded Gaussian mixture model (G2M2). The concentration $c=a+b$ sets how graded the population is. As $c\to0$ the weights pile up at $0$ and $1$, and the law of $X$ becomes a two-component Gaussian mixture with hard memberships. As $c$ grows, the memberships concentrate around their mean and every unit becomes nearly the same blend. When the types' draws are independent, any one-dimensional linear projection of the multivariate models cited above has exactly the form~\eqref{eq:model}.

The model is identifiable, unequal within-type variances included (Section~\ref{sec:ident}). The proof, however, reads the parameters off the moment generating function $M(z)=\E e^{zX}$ as $z\to\infty$, one per order of the expansion, that is, from the far tails of the density, where a sample of moderate size has almost no observations. G2M2 is identifiable and weakly identified at the same time, and most of this paper measures how weak.

\begin{figure*}[t]
\centering
\includegraphics[width=\textwidth]{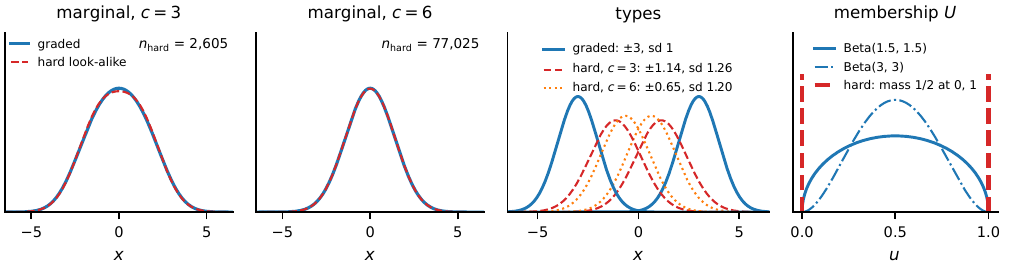}
\caption{Two stories, one marginal ($d=3$). \emph{First two panels:} the graded marginal (solid) and its hard-membership look-alike from Proposition~\ref{prop:lookalike} (dashed) at $c=3$ and $c=6$. The two agree in their first five cumulants, and below the certified sizes shown no test separates them with summed error under $\tfrac12$ (Theorem~\ref{thm:lecam}). \emph{Third panel:} the types behind each story, at $\pm3$ with standard deviation $1$ for the graded model, against $\pm1.14$ with standard deviation $1.26$ ($c=3$) and $\pm0.65$ with standard deviation $1.20$ ($c=6$) for the look-alikes. \emph{Right:} the membership laws, $\Beta(c/2,c/2)$ against mass $\tfrac12$ at each of $0$ and $1$.}
\label{fig:stories}
\end{figure*}

The main result concerns graded versus hard membership (Section~\ref{sec:hard}). For the symmetric model we write down, in closed form, a two-component Gaussian mixture with hard memberships whose first five cumulants equal those of the graded marginal (Figure~\ref{fig:stories}). A Le Cam argument turns the construction into a certified sample size, below which every test of graded against hard membership has type-I and type-II errors summing to at least $\tfrac12$. The certified size depends on $c$ and on $d$, the half-distance between the type means measured in within-type standard deviations, and it changes character at $c_0(d)=4d^4/(4d^2+1)\approx d^2$. At $d=3$ it is 2,605 observations when $c=3$ and 77,025 when $c=6$, and for $c\ge c_0$ it exceeded $6\times10^4$ in every case we computed. Away from $c=0$, where graded and hard membership coincide, the more graded the population, the more data one needs to show that it is graded at all.

\paragraph{Contributions.}
\begin{itemize}\setlength\itemsep{1pt}
\item Identifiability with unequal variances, from a tail expansion that also shows why identification is weak (Section~\ref{sec:ident}).
\item A closed-form hard-membership look-alike matching five cumulants, certified sample sizes, with the sample size needed for reliable testing pinned down to within a factor of 2.6, the regime boundary $c_0(d)$, and a corollary that $c$ cannot be estimated below the certified size (Section~\ref{sec:hard}); and the matching statement at the Gaussian end (Section~\ref{sec:gauss}).
\item Evidence that learning $c$ causes most of the difficulty: it inflates Cram\'er--Rao sample sizes for the separation by a factor of 29 at $c=1$ and by more than $10^4$ for $4.5\le c\le6$ (Section~\ref{sec:fixedc}).
\item Evidence that in practice the prior settles the question: under the default priors of a published sampler for this model class, a correctly specified posterior reports hard membership on all 100 simulated datasets with graded membership, and the sampler's own answers depend on where it starts (Section~\ref{sec:inference}).
\end{itemize}

%% file: sections/model.tex
\section{Model and Related Work}\label{sec:model}

\subsection{The graded Gaussian mixture}

Given $U=u$, the observation in~\eqref{eq:model} is Gaussian,
\begin{equation}\label{eq:conditional}
X\mid U=u\;\sim\;\N\big(u\mu_1+(1-u)\mu_2,\ u^2\sigma_1^2+(1-u)^2\sigma_2^2\big).
\end{equation}
The conditional mean is linear in the membership and the conditional variance is quadratic. The quadratic variance distinguishes averaging independent draws from selecting one of them, as a finite mixture does. We write $c=a+b$ for the concentration and $\pi=a/c$ for the prior mean membership. As $c\to0$ with $\pi$ fixed, $U$ converges in distribution to a Bernoulli$(\pi)$ variable, and the law of $X$ converges to the hard-membership mixture $\pi\,\N(\mu_1,\sigma_1^2)+(1-\pi)\,\N(\mu_2,\sigma_2^2)$. Hard membership is therefore a boundary point of the graded model, not a separate model.

\paragraph{The symmetric model.}
Much of the analysis uses type means $\pm s$, a common within-type standard deviation $t$ and a symmetric membership law, $a=b=r$, so that $c=2r$. Writing $W=2U-1$ and $Z\sim\N(0,1)$ independent of $W$,
\begin{equation}\label{eq:symmetric}
X = sW + t\sqrt{(1+W^2)/2}\;Z .
\end{equation}
We write $\theta=(s,r,t)$, $f_\theta$ for the density of $X$, $v$ for its variance, $\kappa_k$ for its $k$th cumulant and $\gamma_2=\kappa_4/v^2$ for its excess kurtosis. The separation-to-noise ratio is $d=s/t$, so the type means lie $2d$ within-type standard deviations apart. Because $f_{(s,r,t)}(x)=t^{-1}f_{(d,r,1)}(x/t)$, the divergence between $f_\theta$ and any comparison law built from $\theta$ in a scale-equivariant way (its moment-matched Gaussian, or the look-alike of Section~\ref{sec:hard}) depends on $\theta$ only through $(d,c)$. All thresholds below are reported as functions of $(d,c)$.

\subsection{Related work}

\paragraph{Instances of the model.}
\citet{marco2024functional} introduced functional mixed-membership models in which each observation is a convex combination of dependent latent Gaussian features (their term for types), so that the conditional covariance is quadratic in the memberships; \citet{marco2024flexible} treat the high-dimensional multivariate case. Both are implemented in the \texttt{BayesFMMM} package, which we use in Section~\ref{sec:inference}. G2M2 is their model with one-dimensional observations, two features and no cross-covariance between them. They point out an exact rescaling problem: given the memberships, shrinking them into a sub-simplex and extrapolating the feature means leaves the likelihood unchanged. They resolve it by requiring at least one unit to belong entirely to each feature (separability), and their consistency results condition on known memberships. Under a Beta membership law a rescaled membership is no longer Beta distributed, so the exact rescaling is unavailable. The weak direction we study is its closest analogue inside the model: a path along which the low-order moments of the marginal stay nearly fixed, because raising $c$, which compresses memberships toward $\tfrac12$, is offset by moving the type means apart.

In hyperspectral imaging, the normal compositional model treats each pixel as a convex combination of random endmember spectra \citep{stein2003,eismann2007}, with Bayesian estimation developed by \citet{eches2010}, \citet{halimi2015} and \citet{zhou2018}. \citet{woodbridge2019} consider estimating $K$ Gaussians from convex combinations of their realizations, which is our setting with $K$ types. They note that Bayesian estimators for this model are sensitive to their priors, for which our results give a quantitative reason.

Two related models do not contain G2M2: the continuous mixed-membership model of \citet{gruhl2014} draws each outcome from a single type, and the partial-membership model of \citet{heller2008} gives a precision-weighted conditional law.

\paragraph{Identifiability and resolution.}
Finite mixtures of Gaussians are identifiable up to relabeling \citep{teicher1963,yakowitz1968} but can be weakly identified: \citet{ho2016} show that location--scale Gaussian mixtures are estimated slowly because derivatives in the mean and in the variance are linked by the heat equation, and \citet{nguyen2013} and \citet{heinrich2018} relate estimation rates to the strength of identifiability. Testing the number of components is non-regular because the null lies on the boundary of the parameter space \citep{ghosh1985,chen2001}. Hard membership is such a boundary point of G2M2, and a singularity in the sense of \citet{watanabe2009}. \citet{andrews2012} treat inference when the strength of identification varies over the parameter space. For mixed membership, \citet{yurochkin2019} give geometric estimators for the Dirichlet simplex nest, in which noise is additive and does not scale with the membership, and \citet{gu2023} establish identifiability of grade-of-membership models for multivariate categorical data. Neither covers draw averaging. Identifiability does not settle how much data a problem needs. \citet{hardt2015} show that estimating the parameters of a mixture of two univariate Gaussians to constant additive error takes $\Theta(\sigma^{12})$ samples, where $\sigma^2$ is the mixture variance; their lower bound adds Gaussian noise to two mixtures that agree in their first five moments. \citet{wu2020} use moment matching to prove optimal rates for Gaussian mixtures. Our look-alike is a construction of the same kind, but across two model classes: a graded marginal and a hard mixture that agree through the fifth cumulant. The lower bounds use the two-point method of Le Cam \citep{lecam1973,tsybakov2009}.

%% file: sections/identifiability.tex
\section{Identifiability}\label{sec:ident}

The moment generating function $M(z)=\E e^{zX}$ is finite on $\R$ because $U$ is bounded. Conditioning on $U$ in~\eqref{eq:conditional} writes $M(z)$ as a Beta integral of $\exp\{\varphi_z(u)\}$, with $\varphi_z$ quadratic in $u$. For large $z$ the integral concentrates at the endpoint of $[0,1]$ attached to the type with the larger variance (or, if the variances are equal, the larger mean), and a Laplace-type expansion at that endpoint gives the following.

\begin{lemma}[Tail of the moment generating function]\label{lem:mgf2}
In model~\eqref{eq:model} with $(\mu_1,\sigma_1)\ne(\mu_2,\sigma_2)$, label the types so that $(\sigma_1,\mu_1)$ exceeds $(\sigma_2,\mu_2)$ lexicographically, and let $\Delta=\mu_1-\mu_2$. As $z\to+\infty$,
\[
\log M(z)=\tfrac12\sigma_1^2z^2+\mu_1z-2b\log z+C_0+\frac{C_1}{z}+\frac{C_2}{z^2}+O(z^{-3}),
\]
where $C_0=\log\{\Gamma(b)/B(a,b)\}-b\log\sigma_1^2$, $C_1=-b\Delta/\sigma_1^2$, and
\[
C_2=\frac{b\Delta^2}{2\sigma_1^4}+\frac{b(b+1)(\sigma_1^2+\sigma_2^2)}{2\sigma_1^4}-\frac{(a-1)b}{\sigma_1^2}.
\]
\end{lemma}

The proof (Appendix~\ref{app:ident}) substitutes $u=1-\tau/B_z$ with $B_z=\sigma_1^2z^2+\Delta z$. This rescaling absorbs the $z^2$ term of $\varphi_z$ into the Gamma kernel $e^{-\tau}\tau^{b-1}$, and a second-order expansion in the manner of Watson's lemma follows. The coefficients were also checked numerically at 50-digit precision (Appendix~\ref{app:ident}).

\begin{proposition}[Identifiability]\label{prop:ident}
Among parameter vectors of model~\eqref{eq:model} with $(\mu_1,\sigma_1)\ne(\mu_2,\sigma_2)$, the law of $X$ determines $(\mu_1,\mu_2,\sigma_1,\sigma_2,a,b)$ up to the label swap $(\mu_1,\mu_2,\sigma_1,\sigma_2,a,b)\mapsto(\mu_2,\mu_1,\sigma_2,\sigma_1,b,a)$. In particular $c$ is identified, and $\pi$ is identified up to $\pi\mapsto1-\pi$.
\end{proposition}

\begin{proof}[Proof sketch]
The label swap preserves the law because $1-U\sim\Beta(b,a)$, and exactly one labeling satisfies the ordering in Lemma~\ref{lem:mgf2}. Equal laws give equal coefficients in the asymptotic scale $z^2,z,\log z,1,z^{-1},z^{-2}$, which determine in turn $\sigma_1$, $\mu_1$, $b$, then $a$ (since $a\mapsto\log B(a,b)$ is strictly decreasing), $\Delta$, and finally $\sigma_2^2$, which enters $C_2$ linearly with positive coefficient.
\end{proof}

The excluded case $(\mu_1,\sigma_1)=(\mu_2,\sigma_2)$ is a pure scale mixture. Proposition~\ref{prop:ident} does not rule out that a law in its class also arises from the excluded case. Within the symmetric model~\eqref{eq:symmetric}, a separate one-term expansion covers the excluded case as well: $f_\theta=f_{\theta'}$ if and only if $(|s|,r,t)=(|s'|,r',t')$, for all $s\in\R$ (Appendix~\ref{app:ident}).

\begin{remark}[An identification ladder]\label{rem:ladder}
The parameters enter the tail of $\log M$ at successively weaker orders, from $\sigma_1$ at $z^2$ down to the smaller within-type variance at $z^{-2}$, and every step takes $z\to\infty$. For roughly Gaussian data with variance $v$ the empirical moment generating function has relative variance of order $e^{vz^2}/n$, so a sample of size $n$ sees $M(z)$ only for $z$ up to about $\sqrt{(\log n)/v}$, far from where the expansion separates the parameters.
\end{remark}

%% file: sections/graded_vs_hard.tex
\section{Graded Versus Hard Membership}\label{sec:hard}

This section works in the symmetric model~\eqref{eq:symmetric}. Proofs are in Appendix~\ref{app:hard}.

\subsection{Cumulants and the sign of the kurtosis}

\begin{lemma}[Cumulants]\label{lem:cumulants}
All odd cumulants of $f_\theta$ vanish, and
\[
v=\frac{s^2+(r+1)t^2}{2r+1},\qquad
\kappa_4=\frac{-6s^4+12rs^2t^2+3rt^4}{(2r+1)^2(2r+3)} .
\]
\end{lemma}

The variance is constant along the curves $s^2+(r+1)t^2=(2r+1)v$, so the second moment alone cannot separate separation, noise and concentration. The fourth cumulant can, and its sign decides what a hard-membership impostor must look like.

\begin{proposition}[Range and sign of the excess kurtosis]\label{prop:range}
For every $\theta$, $-2<\gamma_2<3$. Moreover $\gamma_2<0$ if and only if $c<c_0(d)$, where
\[
c_0(d)=\frac{4d^4}{4d^2+1}=d^2-\tfrac14+O(d^{-2}),
\]
and $\gamma_2=0$ at $c=c_0(d)$.
\end{proposition}

Below $c_0$ the marginal is platykurtic, like a two-component mixture with separated means. Above it the marginal is leptokurtic, like a scale mixture. At $d=3$, $c_0=8.76$. Numerically $\sup_\theta\gamma_2\approx0.24$, so the upper bound of $3$ is far from attained; it is the bound that keeps the look-alike below well defined.

\subsection{A hard-membership look-alike}

\begin{proposition}[Look-alike]\label{prop:lookalike}
Let $v$ and $\gamma_2$ be the variance and excess kurtosis of $f_\theta$. Define a two-component Gaussian mixture $G_\theta$ as follows. If $\gamma_2<0$,
\[
G_\theta=\tfrac12\N(s_G,t_G^2)+\tfrac12\N(-s_G,t_G^2),
\]
with $s_G^2=v\sqrt{-\gamma_2/2}$ and $t_G^2=v-s_G^2$. If $\gamma_2=0$, $G_\theta=\N(0,v)$. If $\gamma_2>0$,
\[
G_\theta=\tfrac12\N\big(0,v_+\big)+\tfrac12\N\big(0,v_-\big),\quad v_\pm=v\big(1\pm\sqrt{\gamma_2/3}\big).
\]
Then $G_\theta$ is well defined for every $\theta$, and $\kappa_k(G_\theta)=\kappa_k(f_\theta)$ for $k=1,\dots,5$. For $\gamma_2<0$, $G_\theta$ is the hard-membership limit ($c\to0$) of the symmetric model with half-separation $s_G$ and noise $t_G$: every unit belongs to exactly one type.
\end{proposition}

The match is less remarkable than it may sound: both laws are symmetric, so odd cumulants vanish, and a symmetric two-component mixture has exactly two free parameters with which to match $\kappa_2$ and $\kappa_4$. What matters is that this is enough: the moment-matched Gaussian first differs from $f_\theta$ in $\kappa_4$, and $G_\theta$ only in $\kappa_6$. All numbers below use exact divergences.

\subsection{Certified sample sizes}

\begin{theorem}[Graded against hard]\label{thm:lecam}
Let $\nhard(\theta)=1/\{2\KL(f_\theta\|G_\theta)\}$, which is positive and finite and depends on $\theta$ only through $(d,c)$. For every $n$ and every test $\phi$ based on $n$ i.i.d.\ observations,
\[
f_\theta^{\otimes n}(\phi=1)+G_\theta^{\otimes n}(\phi=0)\;\ge\;1-\tfrac12\sqrt{n/\nhard(\theta)} .
\]
In particular, for $n\le\nhard(\theta)$ the two error probabilities sum to at least $\tfrac12$.
\end{theorem}

The proof is Pinsker's inequality for the product measures, and it holds for any fixed comparison law. A \emph{certified sample size} is a value of $1/(2\KL)$ for an explicit comparison law; we write $\ngauss(\theta)$ for the one against $\N(0,v)$, and use numerically optimized two-component mixtures where they are closer than the closed form.

\paragraph{Tightness.}
For the pair $(f_\theta,G_\theta)$, with squared Hellinger distance $H^2=\tfrac12\int(\sqrt{f_\theta}-\sqrt{G_\theta})^2$, every test has summed error at least $\tfrac12(1-H^2)^{2n}$ and the likelihood-ratio test at most $(1-H^2)^n$ (Appendix~\ref{app:hard}), so the sample size at which summed error $0.05$ becomes achievable is pinned down to within a factor of 2.6. At $d=3$, $H^2/\KL\in[0.250,0.252]$ for $c\in\{1,2,3,4,6,8\}$, which places it between about $9\,\nhard$ and $24\,\nhard$. Against the composite family of all hard mixtures testing can only be harder, so the certified size remains a lower bound there.

\paragraph{Estimating the concentration.}
The testing bound transfers to estimation, because the hard look-alike is a limit of graded models with vanishing concentration.

\begin{corollary}\label{cor:estimation}
Let $\gamma_2(\theta)<0$ and $n\le\nhard(\theta)$. For every estimator $\hat c$ and every $\eta>0$ there is $\epsilon>0$ such that
\[
\max_{\theta'\in\{\theta,\,(s_G,\epsilon/2,t_G)\}}\mathbb P_{\theta'}\big(|\hat c-c(\theta')|\ge(c-\epsilon)/2\big)\ge\tfrac14-\eta .
\]
\end{corollary}

No procedure can tell concentration $c$ from concentration near zero with probability better than $\tfrac34$ under both.

\begin{table}[t]
\caption{Certified sample sizes $1/(2\KL)$ at $d=3$. The look-alike column gives the means and, in parentheses, the standard deviations of $G_\theta$; at $c=12>c_0(3)$ it is the equal-mean scale mixture. ``Optimized'' uses a numerically optimized two-component mixture and is never smaller than the closed form. $\ngauss$ is against $\N(0,v)$.}
\label{tab:d3}
\centering\small\setlength{\tabcolsep}{3.5pt}
\begin{tabular}{rrlrrr}
\toprule
 & & & \multicolumn{2}{c}{$\nhard$} & \\
\cmidrule(lr){4-5}
$c$ & $\gamma_2$ & look-alike $G_\theta$ & closed & optimized & $\ngauss$ \\
\midrule
1  & $-0.976$ & $\pm1.92\ (1.26)$ & 136 & 195 & 9 \\
2  & $-0.620$ & $\pm1.43\ (1.27)$ & 647 & 869 & 33 \\
3  & $-0.403$ & $\pm1.14\ (1.26)$ & 2{,}605 & 3{,}285 & 97 \\
4  & $-0.262$ & $\pm0.93\ (1.24)$ & 9{,}167 & 10{,}959 & 266 \\
6  & $-0.101$ & $\pm0.65\ (1.20)$ & 77{,}025 & 84{,}605 & 2{,}139 \\
8  & $-0.019$ & $\pm0.39\ (1.18)$ & 235{,}951 & 241{,}978 & 53{,}363 \\
12 & $+0.047$ & $0\ (1.18,\,1.04)$ & 187{,}522 & $1.5\times10^{8}$ & 10{,}421 \\
\bottomrule
\end{tabular}
\end{table}

Table~\ref{tab:d3} gives the numbers at $d=3$. For $c\le6$, the certified size against hard membership is 15 to 36 times that against a single Gaussian. Near $c_0$ the Gaussian itself becomes a good impostor, because $\gamma_2$ passes through zero, and $\ngauss$ jumps to 53,363 at $c=8$. Section~\ref{sec:gauss} returns to that end of the model.

\subsection{The regime map}

\begin{figure*}[t]
\centering
\includegraphics[width=\textwidth]{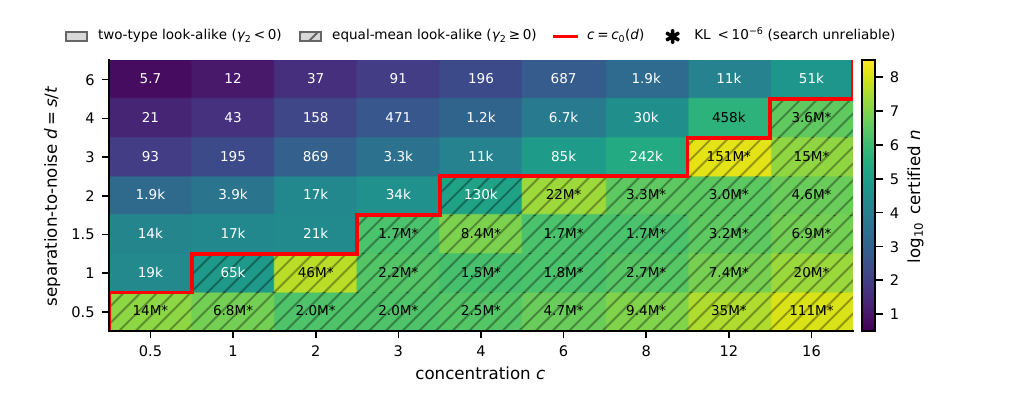}
\caption{Certified sample size below which no test separates graded from hard membership with summed error under $\tfrac12$, across separation-to-noise ratio $d=s/t$ and concentration $c$. Each cell reports $1/(2\KL(f_\theta\|G))$ for the closer, in divergence, of the closed-form look-alike and a numerically optimized two-component mixture. Red staircase: $c=c_0(d)=4d^4/(4d^2+1)\approx d^2$. In the hatched cells to its right, $\gamma_2\ge0$ and the look-alike is an equal-mean scale mixture, so the comparison is with one type of heterogeneous variance rather than two separated types. Every computed cell with $c\ge c_0(d)$ exceeds $6\times10^4$. KL values carry quadrature error below $10^{-15}$. In cells marked $*$ the divergence is below $10^{-6}$ and the numerical search for the closest mixture is unreliable: each value is still a valid lower bound, but these cells need not be monotone in $c$ and should not be compared with one another.}
\label{fig:regime}
\end{figure*}

Figure~\ref{fig:regime} extends the table to a grid over $(d,c)$. The certified size is U-shaped in $c$. As $c\to0$ it diverges, because the graded model itself converges to a hard mixture; this end is benign, since there graded membership is nearly hard. Away from it the certified size rises steeply with $c$; for $d\ge1$ the minimum lies at or below the smallest computed concentration, $c=0.5$. The boundary $c_0(d)\approx d^2$ separates two regimes. To its left, well-separated types at modest concentration have small certified sizes: at $d=6$, types twelve within-type standard deviations apart, the certified size at $c=6$ is 687. To its right, every computed cell exceeds $6\times10^4$, and most exceed $10^6$. When the types are close, the question is hard throughout: for $d\le1.5$ no computed cell falls below $1.4\times10^4$.

\paragraph{Asymmetric truths.}
Theorem~\ref{thm:lecam} applies to any truth and any comparison law. For three asymmetric truths at $d=3$ ($\pi=0.3$, $\sigma_1=1.5$, or both), certified sizes against numerically optimized hard mixtures range from 1,009 to 4,312 at $c=3$ and from 4,462 to 34,218 at $c=6$, about 40 to 125 times the corresponding Gaussian thresholds (Table~\ref{tab:asym}).

%% file: sections/gaussian_end.tex
\section{The Gaussian End}\label{sec:gauss}

Section~\ref{sec:hard} concerned the end of the model where $c\to0$. At the other end, a large concentration pulls every membership toward $\tfrac12$ and the marginal approaches a Gaussian. By Lemma~\ref{lem:cumulants}, the set of symmetric configurations with a given variance $v$ and noise level $t$ contains the point $(s_r,r,t)$ for every $r$ at which
\[
s_r^2=r(2v-t^2)+(v-t^2)\ \ge0 .
\]
When $t^2<2v$ this curve is unbounded, and $s_r\to\infty$ as $r\to\infty$. Along it the marginal converges to $\N(0,v)$.

\begin{proposition}[The Gaussian end]\label{prop:gauss-end}
Fix $\theta$ with variance $v$ and $t^2<2v$, and let $\theta_k=(s_{r_k},r_k,t)$ with $r_k\to\infty$. Then $f_{\theta_k}\to\N(0,v)$ in total variation, and for every $n$ and every test $\phi$,
\[
\liminf_{k\to\infty}\big[f_\theta^{\otimes n}(\phi=1)+f_{\theta_k}^{\otimes n}(\phi=0)\big]\ge1-\tfrac12\sqrt{n/\ngauss(\theta)} .
\]
\end{proposition}

The proof (Appendix~\ref{app:hard}) combines a Beta central limit theorem, Scheff\'e's lemma and Pinsker's inequality. For $n\le\ngauss$, no test can exclude the unbounded part of the curve with summed error below $\tfrac12$. How fast the marginal approaches the Gaussian depends on the path.

\begin{proposition}[Two regimes for the kurtosis]\label{prop:scaling}
(a) Along the curve above, with $\rho=t^2/v\in(0,2)$ fixed, $\gamma_2=3(2-\rho)(3\rho-2)/(4r)+O(r^{-2})$.
(b) At fixed $(s,t)$, $\gamma_2=(6d^2+3/2)/r^2+O(r^{-3})$ as $r\to\infty$.
\end{proposition}

With the heuristic $\KL(f_\theta\|\N(0,v))\approx\gamma_2^2/48$, accurate to 4\% for $r\ge5$ along such curves, $\ngauss$ grows like $c^2$ along a curve of fixed variance (unless $\rho=\tfrac23$) and like $c^4$ at fixed $(s,t)$, where it is also non-monotone, peaking near $c_0(d)$ (Table~\ref{tab:d3}). Neither power law holds over practical ranges of $c$.

Throughout Table~\ref{tab:d3} the hard end binds ($\nhard>\ngauss$). The Gaussian end matters in practice because a sampler started far from the truth tends to settle there (Section~\ref{sec:inference}).

%% file: sections/fixed_c.tex
\section{Supplying the Concentration}\label{sec:fixedc}

\begin{figure*}[t]
\centering
\includegraphics[width=\textwidth]{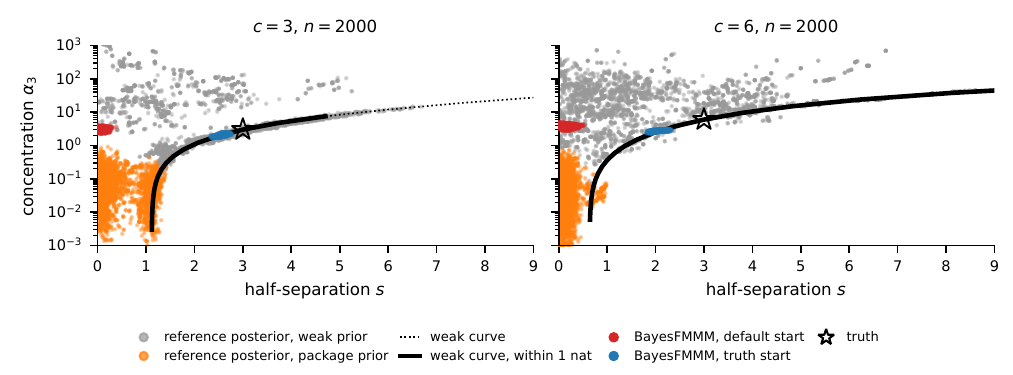}
\caption{Where the posterior puts its mass (dataset 1 at each $c$; $n=2000$, $d=3$). Gray: reference posterior with memberships integrated out, under a weak prior. Orange: the same under \texttt{BayesFMMM}'s priors on the means and on the concentration, $\alpha_3\sim\mathrm{Exp}(10)$. Each pools two runs, started at the truth and at the configuration in which \texttt{BayesFMMM}'s default start settles; at $c=3$ the gray cloud at small $s$ and large $\alpha_3$ comes from the collapsed start, which mixes incompletely. Dotted: the weak curve, the symmetric configuration closest in divergence to the truth at each concentration; thick where it lies within one nat of the truth. Red and blue: \texttt{BayesFMMM} draws from the second half of its long runs, from its default and truth-initialized starts. Star: truth.}
\label{fig:posterior}
\end{figure*}

In applications $c$ is usually learned along with everything else. To see how much of the difficulty that causes, we computed Cram\'er--Rao bounds for the full two-type model at the symmetric truth ($\mu=\pm3$, $\sigma_1=\sigma_2=1$, $\pi=\tfrac12$), with $\mu_1,\mu_2,\sigma_1,\sigma_2$ and $\pi$ free and $c$ either supplied or learned. Table~\ref{tab:fisher} reports, for each target, the sample size at which the bound on the standard deviation equals a stated tolerance.

\begin{table}[t]
\caption{Cram\'er--Rao sample sizes at $d=3$: the $n$ at which the bound on the standard deviation equals the tolerance. Tolerances are $0.3$ (10\%) for the half-separation, 10\% for $\sigma_1$, $0.05$ for $\pi$ and $0.1$ for $\log c$. The last column is the ratio of the two half-separation columns, computed from unrounded bounds.}
\label{tab:fisher}
\centering\small\setlength{\tabcolsep}{3pt}
\begin{tabular}{r|rrr|rr|r}
\toprule
 & \multicolumn{3}{c|}{$c$ supplied} & \multicolumn{2}{c|}{$c$ learned} & \\
$c$ & half-sep. & $\sigma_1$ & $\pi$ & half-sep. & $\log c$ & ratio \\
\midrule
1 & 44 & 754 & 220 & 1{,}259 & 13{,}266 & 29 \\
2 & 80 & 2{,}144 & 1{,}077 & 12{,}605 & 75{,}167 & 158 \\
3 & 132 & 5{,}234 & 5{,}266 & 95{,}022 & 458{,}203 & 718 \\
4 & 218 & 13{,}017 & 25{,}634 & 810{,}078 & $3.5\!\times\!10^{6}$ & 3{,}719 \\
6 & 576 & 225{,}012 & $1.6\!\times\!10^{6}$ & $9.9\!\times\!10^{6}$ & $4.1\!\times\!10^{7}$ & 17{,}236 \\
8 & 1{,}408 & 36{,}407 & 805{,}678 & $2.8\!\times\!10^{6}$ & $1.2\!\times\!10^{7}$ & 1{,}969 \\
\bottomrule
\end{tabular}
\end{table}

With $c$ supplied, the separation is cheap: 44 to 1,408 observations bring the bound on its standard deviation down to 10\%. Learning $c$ multiplies these sizes by 29 at $c=1$ and by about 17,000 at $c=6$; on a finer grid the factor peaks at about 54,000 near $c=5.5$ (Appendix~\ref{app:numerics}). Supplying $c$ removes the curve of Section~\ref{sec:gauss}, since at fixed $r$ and $v$ the relation $s^2+(r+1)t^2=(2r+1)v$ bounds $s$.

Something weak remains. At $c=6$ the weakest Fisher direction, in coordinates $(\mu_1,\mu_2,\log\sigma_1,\log\sigma_2,\operatorname{logit}\pi)$, is about $(0.64,0.64,-0.08,0.08,-0.42)$: it shifts both means together, moves the log-variances in opposite directions and changes $\pi$, leaving the separation almost unchanged. This is why $\sigma_1$ and $\pi$ are expensive at $c=6$ (225,012 and $1.6\times10^6$ observations) although the separation is not. These columns are non-monotone in $c$; on a finer grid both peak near $c=6.5$.

\paragraph{A check against fitted estimates.}
We fitted the full model, with the true $c$ supplied, to ten datasets per cell at $n\in\{10^3,10^4\}$ by maximizing the marginal likelihood from many starts (Appendices~\ref{app:sweep} and~\ref{app:em}). At $n=10^4$ the across-dataset standard deviation of $\hat s$ is within the sampling band of the bound for $0.25\le c\le4$ (ratios 0.71 to 1.27), but about 2.6 times the bound at $c=6$ and $8$, where a few fits sit on other local maxima. For $c\ge10$ the likelihood has competing maxima within about a nat of one another, collapsed, separated and point-mass configurations, and the fit depends on the number of starts. At $n=10^3$, four to eight of ten datasets had a rival maximum within one nat whose half-separation differed by more than 0.5. There the local bound no longer describes the estimator; this is the weak identification of Sections~\ref{sec:hard} and~\ref{sec:gauss} seen from the estimator's side.

%% file: sections/inference.tex
\section{Existing Inference}\label{sec:inference}

We studied the published sampler for this model class, \texttt{BayesFMMM} \citep{marco2024functional,marco2024flexible}, on data from the symmetric model with $d=3$, $n=2000$ and $c\in\{3,6\}$, ten datasets per value of $c$. With one-dimensional observations its noise variance is absorbed into the feature covariance, so G2M2 lies inside the model it fits, although scalar data are far from the functional and high-dimensional settings it was built for. We ask two questions: what a correct posterior under the package's default priors reports, and whether the package's sampler finds it. Configurations are compared by the exact marginal log-likelihood, not the conditional one the sampler prints.

\paragraph{The prior decides.}
To answer the first question we sampled the same model with memberships integrated out by quadrature, using an ensemble sampler \citep{foreman2013}. Each dataset was run from the truth and from the configuration where \texttt{BayesFMMM}'s default start settles, under two priors: the package's priors on the means and on $\alpha_3$, which is exponential with mean $0.1$, and a weak prior (Appendix~\ref{app:bayes}). On 100 datasets, 50 per value of $c$, the package's priors put the posterior at the hard end every time and from both starts. The pooled median of $\alpha_3$ is 0.03 to 0.09, every 95\% interval excludes the true concentration, and the half-separation is between about 0.1 and 1.2, against 3 in truth. Yet the best hard configuration fits about as well as the truth: the median log-likelihood difference is between $-0.7$ and $+2.0$ nats across values of $c$ and starts. Under the weak prior no dataset goes hard and the intervals cover the truth (all for $\alpha_3$, 192 of 200 for $s$), but their upper limits for $\alpha_3$ range from 5.5 to 907. At $n=2000$ the data do not determine whether membership is graded; the priors do (Figure~\ref{fig:posterior}). With the package's priors this persists far beyond the certified sizes: at $n=120{,}000$ for $c=3$ and $n=10^6$ for $c=6$ the intervals still exclude the true concentration and separation. A rough accounting attributes this to the priors placing the truth tens of nats below a hard configuration, while the data's evidence grows only like $n\KL$ (Appendix~\ref{app:bayes}).

\paragraph{The sampler.}
We then ran \texttt{BayesFMMM} itself under its default priors, from its default initialization and from the truth, in short runs (2,000 iterations, one initialization try, ten datasets per $c$) and in long runs on five of them (50,000 iterations, multiple starts, tempered transitions; Appendix~\ref{app:bayes}).

\begin{table}[t]
\caption{\texttt{BayesFMMM} at $d=3$, $n=2000$, default priors, two starting points. Short runs: 2,000 iterations, one initialization try, ten datasets per $c$. Long runs: 50,000 iterations, multiple starts, tempered transitions, five of the same datasets. ``Gap'' is the marginal log-likelihood of the truth minus that of the best of 12 posterior draws, memberships integrated out; the predicted gap to the Gaussian end, $n\KL(f_\theta\|\N(0,v))$, is 10.3 at $c=3$ and 0.47 at $c=6$. Coverage is of the nominal 95\% interval.}
\label{tab:rep}
\centering\small\setlength{\tabcolsep}{3pt}
\begin{tabular}{lrrrr}
\toprule
 & \multicolumn{2}{c}{$c=3$} & \multicolumn{2}{c}{$c=6$} \\
\cmidrule(lr){2-3}\cmidrule(lr){4-5}
start & default & truth & default & truth \\
\midrule
\multicolumn{5}{l}{\emph{short runs}} \\
median $\hat s$ (truth 3) & 0.10 & 2.97 & 0.15 & 2.95 \\
median $\hat\alpha_3$ (truth $c$) & 12.7 & 2.9 & 12.7 & 5.5 \\
covers $s$, $\alpha_3$ (\%) & 0, 0 & 80, 60 & 0, 0 & 70, 20 \\
gap, median & $+8.9$ & $-1.1$ & $+0.4$ & $-1.0$ \\
\multicolumn{5}{l}{\emph{long runs}} \\
median $\hat s$ & 0.15 & 2.55 & 0.14 & 2.07 \\
median $\hat\alpha_3$ & 3.5 & 2.1 & 3.5 & 2.8 \\
covers $s$, $\alpha_3$ (\%) & 0, 80 & 0, 0 & 0, 0 & 0, 0 \\
gap, median & $+10.9$ & $-1.3$ & $+0.4$ & $-0.9$ \\
\bottomrule
\end{tabular}
\end{table}

From its default start the sampler settles near zero separation, with memberships uncorrelated with the true ones; longer runs change its concentration but not its separation (Table~\ref{tab:rep}). The two starts gave disjoint 95\% intervals for $s$ on all 20 short-run and all 10 long-run datasets (split-$\hat R$ 12 to 26). Tempered transitions hot enough to cross between the configurations were essentially never accepted with 2,000 latent memberships. Nor is the truth start a calibrated reference: over 50,000 iterations it drifts toward smaller $s$ and $\alpha_3$ along the weak curve, the symmetric configurations closest in divergence to the truth at each concentration (Figure~\ref{fig:posterior}), and its intervals stop covering the truth.

The likelihood gaps separate sampler failure from genuine ambiguity. The default-start configuration is close to $\N(0,v)$, so its gap to the truth should be near $n\KL(f_\theta\|\N(0,v))$. That is 10.3 nats at $c=3$, where the observed medians are 8.9 and 10.9 and every gap is positive, so the sampler fails to reach a better fit. At $c=6$ it is 0.47 and the observed median is 0.4, since $n=2000$ is close to $\ngauss=2{,}139$: the data cannot choose.

%% file: sections/discussion.tex
\section{Discussion and Limitations}\label{sec:limits}

G2M2 is identifiable, yet over much of its parameter space the data cannot tell graded from hard membership at the sample sizes one usually has. A practical consequence is a diagnostic. Before reporting that units are blends, compute the certified size against the best-fitting hard mixture over the configurations the data cannot exclude, such as posterior draws or fits within a nat of the best; a point estimate is unreliable in exactly this regime (Section~\ref{sec:fixedc}). If $n$ is below it for any of them, the data cannot rule out hard membership, and a conclusion of gradedness rests on the prior or the sampler's starting point.

\paragraph{Limitations.}
Our results are for scalar data and two types. When the types' draws are independent, any one-dimensional linear projection of a multivariate model of this kind is exactly a G2M2, so the thresholds bound what any single projection can show. With $P$ features sharing one membership, however, features in which the two types have equal means still reveal each unit's variance factor $U^2+(1-U)^2$, and we expect this to weaken the ambiguity as $P$ grows; we have not quantified how fast. The closed-form thresholds are for the symmetric model; for asymmetric truths the certificates are numerical. The case $(\mu_1,\sigma_1)=(\mu_2,\sigma_2)$ is treated only in its symmetric instance, and for $K\ge3$, scalar data constrain a membership only through two conditional moments, so different memberships can induce the same law, an exact obstruction we do not treat. Finally, Section~\ref{sec:inference} evaluates one sampler, and at large $n$ its collapsed-start reference chains do not mix.

\paragraph{Code availability.}
Code for all computations, and the outputs behind every table and figure, are available from the authors; a documented public repository is in preparation.

%% file: sections/appendix.tex
\input{sections/appendix_proofs}

\section{Numerical Methods and Additional Tables}\label{app:numerics}

\paragraph{Densities and divergences.}
Graded densities are computed by Gauss--Jacobi quadrature in the membership, with 100 to 160 nodes. Divergences are computed by Simpson's rule on grids of 2,401 to 4,001 points, with quadrature error in KL below $10^{-15}$. For the optimized look-alikes we minimized $\KL(f_\theta\|G)$ over two-component Gaussian mixtures $G$, parameterized by the logit of the weight, the two means and the two log standard deviations, with the Nelder--Mead method in SciPy. For the $(d,c)$ grid of Figure~\ref{fig:regime} there were 11 starts (symmetric hard mixtures at $\pm\{0.6,1.0,1.4\}\,s_G$ and eight equal-mean scale mixtures), with tolerances $10^{-8}$ in the parameters and $10^{-15}$ in the objective and at most 3,000 evaluations; the reported value uses whichever of the best search and the closed form has the smaller divergence. For the asymmetric truths there were 24 starts and at most 2,500 evaluations. We did not check convergence flags. This does not affect validity: each reported value is the divergence of an explicit mixture, so a suboptimal search only makes the certificate conservative. Where $\KL$ falls below about $10^{-6}$ the search is less reliable, which explains the non-monotone entries in the $\gamma_2\ge0$ region.

\paragraph{Cram\'er--Rao bounds.}
Table~\ref{tab:fisher} uses the Fisher information of the full two-type model at the symmetric truth. The density is computed with 120 Gauss--Jacobi nodes in the membership on a grid of 2,401 points over $[-12,12]$; its derivatives are central differences with step $10^{-4}$ in the coordinates $(\mu_1,\mu_2,\log\sigma_1,\log\sigma_2,\operatorname{logit}\pi,\log c)$, and the information is integrated by Simpson's rule on the same grid. The non-monotonicity in $c$ was checked against changes in step size, number of quadrature nodes and grid width. At $c=5$ and $c=7$ the $\sigma_1$ column is 40,001 and 177,014, and the $\pi$ column is 149,883 and $2.26\times10^6$. The bounds were also computed on a finer grid ($c=0.1$, $0.25$, then $0.5$ to $16$ in steps of $0.5$, and $18$, $20$, $24$). On the finer grid the fixed-$c$ bounds for $\sigma_1$ and $\pi$ both peak near $c=6.5$, and the learned-$c$ bound for the half-separation peaks near $c=5.5$, where it is about 54,000 times the fixed-$c$ bound, before falling. At $c=0.1$ and $0.25$ the Beta density is singular at both endpoints; Gauss--Jacobi quadrature integrates the singularity exactly, and the bounds change by at most $6\times10^{-4}$ in relative terms between 160 and 400 nodes, across step sizes and with a wider grid. With $c$ learned, the bound for the marginal variance equals $(2+\gamma_2)v^2$, the variance of the sample variance, to four digits at every $c$ on the grid; with $c$ supplied it equals it for $c\ge4$ and is up to about 2\% smaller below. As a check between local and global calculations, the learned-$c$ half-separation sizes are about half of the sample sizes at which the best configuration with a 10\% change in separation reaches a divergence of one nat.

\begin{table}[h]
\caption{Tightness of Theorem~\ref{thm:lecam} at $d=3$ for the closed-form look-alike. The certified size is $\nhard=1/(2\KL)$. The last two columns bracket the smallest $n$ at which some test has summed error at most $0.05$ (Proposition~\ref{prop:tight}): no test achieves it below the first, and the likelihood-ratio test achieves it at the second.}
\label{tab:tight}
\centering\small
\begin{tabular}{rrrrrrr}
\toprule
$c$ & $\KL(f_\theta\|G_\theta)$ & $H^2(f_\theta,G_\theta)$ & $H^2/\KL$ & $\nhard$ & lower & upper \\
\midrule
1 & $3.69\times10^{-3}$ & $9.22\times10^{-4}$ & 0.250 & 136 & 1{,}248 & 3{,}247 \\
2 & $7.73\times10^{-4}$ & $1.95\times10^{-4}$ & 0.252 & 647 & 5{,}908 & 15{,}373 \\
3 & $1.92\times10^{-4}$ & $4.84\times10^{-5}$ & 0.252 & 2{,}605 & 23{,}780 & 61{,}876 \\
4 & $5.45\times10^{-5}$ & $1.38\times10^{-5}$ & 0.252 & 9{,}167 & 83{,}718 & 217{,}839 \\
6 & $6.49\times10^{-6}$ & $1.63\times10^{-6}$ & 0.252 & 77{,}025 & 705{,}076 & 1{,}834{,}649 \\
8 & $2.12\times10^{-6}$ & $5.32\times10^{-7}$ & 0.251 & 235{,}951 & 2{,}163{,}447 & 5{,}629{,}418 \\
\bottomrule
\end{tabular}
\end{table}

\begin{table}[h]
\caption{Certified sizes $1/(2\KL)$ for asymmetric truths with $\mu_1=3$, $\mu_2=-3$, against the numerically optimized two-component mixture (hard) and the moment-matched Gaussian.}
\label{tab:asym}
\centering\small
\begin{tabular}{lrrrr}
\toprule
 & \multicolumn{2}{c}{$c=3$} & \multicolumn{2}{c}{$c=6$} \\
truth & hard & Gaussian & hard & Gaussian \\
\midrule
$\pi=0.5$, $\sigma_1=\sigma_2=1$ (symmetric) & 3{,}285 & 97 & 84{,}605 & 2{,}139 \\
$\pi=0.3$, $\sigma_1=\sigma_2=1$ & 4{,}312 & 38 & 34{,}218 & 273 \\
$\pi=0.5$, $\sigma_1=1.5$, $\sigma_2=1$ & 2{,}706 & 67 & 6{,}988 & 98 \\
$\pi=0.3$, $\sigma_1=1.5$, $\sigma_2=1$ & 1{,}009 & 22 & 4{,}462 & 64 \\
\bottomrule
\end{tabular}
\end{table}

\begin{table}[h]
\caption{The Gaussian end along the fixed-$(s,t)$ family $s=3$, $t=1$. The column $1/\chi^2$ is a chi-squared counterpart of $\ngauss=1/(2\KL(f_\theta\|\N(0,v)))$, which is shown in Table~\ref{tab:d3}. $c=8.76$ is $c_0(3)$, where $\gamma_2=0$.}
\label{tab:gauss}
\centering\small
\begin{tabular}{rrrr}
\toprule
$c$ & $\gamma_2$ & $\chi^2(f_\theta\|\N(0,v))$ & $1/\chi^2$ \\
\midrule
1 & $-0.976$ & $1.06\times10^{-1}$ & 9.4 \\
2 & $-0.620$ & $2.76\times10^{-2}$ & 36 \\
3 & $-0.403$ & $9.28\times10^{-3}$ & 108 \\
4 & $-0.262$ & $3.43\times10^{-3}$ & 291 \\
6 & $-0.101$ & $4.40\times10^{-4}$ & 2{,}271 \\
8 & $-0.020$ & $1.81\times10^{-5}$ & 55{,}259 \\
8.76 & $\approx0$ & $4.48\times10^{-6}$ & 223{,}266 \\
10 & $+0.024$ & $2.97\times10^{-5}$ & 33{,}625 \\
16 & $+0.065$ & $1.81\times10^{-4}$ & 5{,}542 \\
\bottomrule
\end{tabular}
\end{table}

\section{The Fixed-Concentration Sweep}\label{app:sweep}

\paragraph{Protocol.}
For each $c\in\{0.1,0.25,0.5,1,2,3,4,6,8,10,15,20\}$ and $n\in\{10^3,10^4\}$ we drew ten datasets from the symmetric model with $s=3$, $t=1$ (seeds 67--76), standardized them, and fitted the full two-type model, with $\mu_1,\mu_2,\sigma_1,\sigma_2$ and $\pi$ free and $c$ fixed at its true value, as described in Appendix~\ref{app:em}: EM from six starts followed by a quasi-Newton polish of the marginal likelihood from the three best endpoints, and for $c\ge6$ EM from forty starts with the ten best distinct endpoints polished. Refits in fresh processes reproduce every fit exactly. Estimates are transformed back to the original scale. We report $\hat s=|\hat\mu_1-\hat\mu_2|/2$, the number of fits in which one type's standard deviation falls below 0.01 in standardized units (a point-mass type), and the variance $\hat v$ implied by the fitted parameters through
\[
v=\sum_k\sigma_k^2\Big(\bar m_k^2+\frac{\bar m_k(1-\bar m_k)}{c+1}\Big)+\frac{\bar m_1\bar m_2(\mu_1-\mu_2)^2}{c+1},\qquad \bar m=(\pi,1-\pi),
\]
which does not assume that the fit stayed symmetric.

\begin{figure}[t]
\centering
\includegraphics[width=\textwidth]{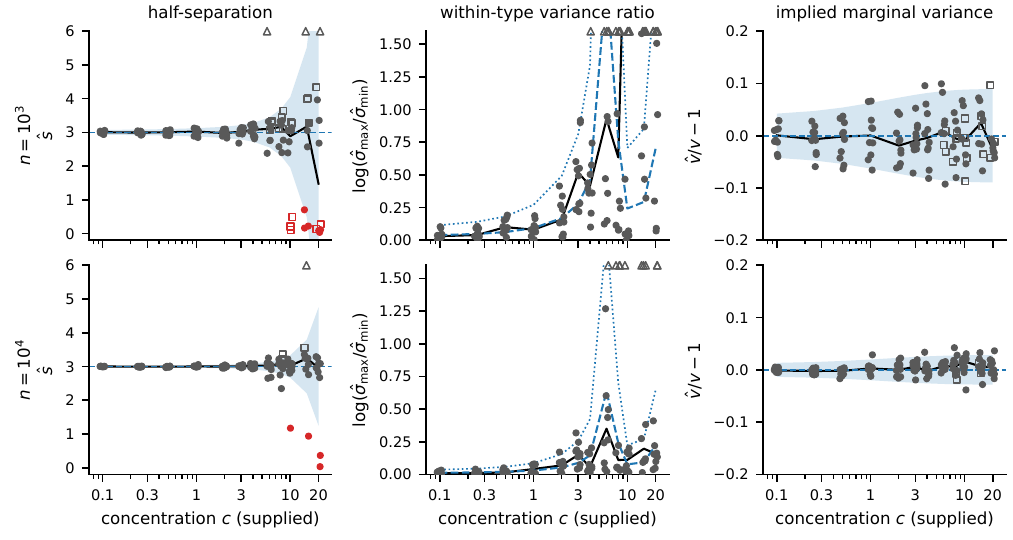}
\caption{Fixed-$c$ sweep against Cram\'er--Rao bounds (ten datasets per cell, truth $s=3$, $t=1$; top row $n=10^3$, bottom row $n=10^4$). Points are individual fits (six starts for $c\le4$, forty for $c\ge6$); red marks $\hat s<1.5$, open squares mark point-mass fits, and open triangles at the top edge mark values beyond the plotted range. Black lines are medians. \emph{Left:} $\hat s$ with the band $3\pm1.96\,\mathrm{sd}$ from the bound. \emph{Middle:} $\log(\hat\sigma_{\max}/\hat\sigma_{\min})$, the absolute value of $\log\hat\sigma_1-\log\hat\sigma_2$, against $0.674\,\mathrm{sd}$ (dashed, the median of a folded normal) and $1.96\,\mathrm{sd}$ (dotted) from the bound. \emph{Right:} $\hat v/v-1$ with the band $\pm1.96\sqrt{(2+\gamma_2)/n}$.}
\label{fig:sweep}
\end{figure}

\begin{table}[h]
\caption{Fixed-$c$ sweep, ten datasets per cell (six starts for $c\le4$, forty for $c\ge6$). sd: across-dataset standard deviation of $\hat s$; CRB: Cram\'er--Rao standard deviation of the half-separation; \#: fits with $\hat s<1.5$; pm: point-mass fits; $\hat v$: relative across-dataset standard deviation of $\hat v$ divided by $\sqrt{(2+\gamma_2)/n}$. $^\dagger$: start-dependent cell, in which forty starts changed several fits relative to six; its values are the best available, not certified maxima. For ten datasets the 95\% sampling band of sd/CRB is $[0.55,1.45]$.}
\label{tab:sweep}
\centering\footnotesize\setlength{\tabcolsep}{3.5pt}
\begin{tabular}{r|rrrrrr|rrrrrr}
\toprule
 & \multicolumn{6}{c|}{$n=10^3$} & \multicolumn{6}{c}{$n=10^4$} \\
$c$ & sd & CRB & med $\hat s$ & \# & pm & $\hat v$ & sd & CRB & med $\hat s$ & \# & pm & $\hat v$ \\
\midrule
0.1 & 0.041 & 0.036 & 3.01 & 0 & 0 & 1.12 & 0.005 & 0.012 & 3.00 & 0 & 0 & 0.52 \\
0.25 & 0.043 & 0.042 & 2.99 & 0 & 0 & 1.10 & 0.009 & 0.013 & 3.00 & 0 & 0 & 0.68 \\
0.5 & 0.064 & 0.050 & 3.00 & 0 & 0 & 1.35 & 0.013 & 0.016 & 2.99 & 0 & 0 & 0.81 \\
1 & 0.091 & 0.063 & 3.02 & 0 & 0 & 1.40 & 0.020 & 0.020 & 3.00 & 0 & 0 & 1.09 \\
2 & 0.084 & 0.085 & 2.98 & 0 & 0 & 1.03 & 0.029 & 0.027 & 3.00 & 0 & 0 & 1.21 \\
3 & 0.131 & 0.109 & 3.02 & 0 & 0 & 0.83 & 0.029 & 0.035 & 3.01 & 0 & 0 & 0.85 \\
4 & 0.175 & 0.140 & 3.08 & 0 & 0 & 1.32 & 0.056 & 0.044 & 3.02 & 0 & 0 & 1.08 \\
6 & 1.377$^\dagger$ & 0.228 & 3.10 & 0 & 3 & 1.25 & 0.187 & 0.072 & 3.03 & 0 & 0 & 0.50 \\
8 & 0.396 & 0.356 & 3.16 & 0 & 4 & 0.85 & 0.287 & 0.113 & 3.08 & 0 & 2 & 1.33 \\
10 & 1.344$^\dagger$ & 0.533 & 2.88 & 3 & 5 & 1.10 & 0.593$^\dagger$ & 0.169 & 3.02 & 1 & 0 & 1.37 \\
15 & 2.804$^\dagger$ & 1.285 & 3.17 & 3 & 3 & 1.00 & 2.393$^\dagger$ & 0.406 & 3.24 & 1 & 2 & 1.18 \\
20 & 2.852$^\dagger$ & 2.825 & 1.48 & 5 & 3 & 0.89 & 1.206$^\dagger$ & 0.893 & 2.96 & 2 & 0 & 1.10 \\
\bottomrule
\end{tabular}
\end{table}

\paragraph{Where the bound describes the fits.}
For $0.25\le c\le4$ the spread of $\hat s$ matches the bound at both sample sizes, with ratios between 0.71 and 1.46. At $c=0.1$, $n=10^4$ it is below the bound (ratio 0.47). This is a property of the ten datasets, not of the estimator: on them the naive estimator (half the difference between the means of the positive and the negative observations) also has standard deviation 0.0063, and on 80 further datasets the EM alone gives 0.0106 against a bound of 0.0115, with blocks of ten datasets ranging from 0.0064 to 0.0137.

\paragraph{Where it does not.}
At $c=6$ and $8$ and $n=10^4$, forty starts reproduce the six-start fits exactly, and the spread is 2.59 and 2.55 times the bound (1.45 and 1.51 for the robust, interquartile estimate): a few fits sit on other local maxima. For $c\ge10$, and at $c=6$ with $n=10^3$, the fits depend on the number of starts. Relative to six starts, forty changed five to seven of ten fits at $n=10^3$ and two of ten at $n=10^4$, by up to 2.0 nats. The competing maxima are qualitatively different configurations: collapsed ($\hat s\le0.5$), separated ($\hat s\approx3$) and point-mass ($\hat s$ of 7 to 10 with one standard deviation near zero). The fraction of datasets with a rival maximum within one nat of the kept one, whose $\hat s$ differs from it by more than 0.5, is

\begin{center}\footnotesize
\begin{tabular}{rccccc}
\toprule
$n$ & $c=6$ & $c=8$ & $c=10$ & $c=15$ & $c=20$ \\
\midrule
$10^3$ & 3/10 & 3/10 & 6/10 & 4/10 & 8/10 \\
$10^4$ & 0/10 & 0/10 & 1/10 & 2/10 & 1/10 \\
\bottomrule
\end{tabular}
\end{center}

\noindent At $c=20$ and $n=10^3$, the number of collapsed fits grows with the number of starts: none with six, three with twenty and five with forty. The collapsed and separated maxima differ by less than a nat. In 14 of the 240 six-start fits one type's standard deviation drifts toward zero. In these fits the Beta exponent at that type's pure end exceeds one, so the likelihood stays bounded and becomes flat in that direction (Appendix~\ref{app:em-code}). Such point-mass configurations belong to the weakly identified set, and some are the highest maximum found, by 1.06 and 1.56 nats in two datasets.

\paragraph{Variances.}
For $c\le4$ the median of $|\log\hat\sigma_1-\log\hat\sigma_2|$ is between 0.39 and 1.85 times the folded-normal prediction $0.674\,\mathrm{sd}$. At $c=6$ and $8$, excluding point-mass fits, it is 0.23 to 0.55 times the prediction. That shortfall is what one expects of a local bound where the likelihood is far from quadratic. Along the weakest Fisher direction at $c=6$, at the displacement where the quadratic approximation reaches one nat, the exact divergence is $2.7\times10^3$ times larger at $n=10^3$ and $9.2\times10^3$ times larger at $n=10^4$. The exact one-nat displacement corresponds to $|\Delta\log(\sigma_1/\sigma_2)|$ of about 0.42 and 0.27, against 4.2 and 1.3 from the quadratic approximation. This is a line search along one direction, not a full profile. The implied marginal variance is estimated as well as a second moment can be in every cell: its relative spread is between 0.50 and 1.40 times $\sqrt{(2+\gamma_2)/n}$.

\section{Details of the Inference Experiments}\label{app:bayes}

\paragraph{\texttt{BayesFMMM}.}
We used the multivariate pipeline of the R package (\texttt{BMVMMM\_Nu\_Z\_multiple\_try}, then \texttt{BMVMMM\_Theta\_est}, then \texttt{BMVMMM\_warm\_start}) with one-dimensional observations ($P=1$), two features ($K=2$) and $M=2=KP$ eigencomponents, so the covariance representation is exact. The default priors are $Z_i\sim\mathrm{Dir}(\alpha_3\pi)$ for the memberships, $\alpha_3\sim\mathrm{Exp}(\text{rate }10)$, $\pi\sim\mathrm{Dir}(\mathbf 1)$, an inverse-gamma$(10,1)$ prior on the variance of the feature means (so the means are marginally $t_{20}$ with scale 0.316), and $\sigma^2\sim\text{inverse-gamma}(1,1)$ for the noise. With one-dimensional observations, $\sigma^2$ is exactly redundant with the feature covariance: since the memberships sum to one, $\sigma^2=\sigma^2(Z_{i1}+Z_{i2})^2$, so only $C+\sigma^2\mathbf 1\mathbf 1^\top$ is identified, and G2M2 corresponds to this matrix being diagonal. The truth-initialized runs overwrite the warm-start objects with the true values. The package's membership rescaling step, which enforces the separability condition of \citet{marco2024functional}, is not exported, and we did not use it.

\emph{Short runs} (ten datasets per $c$, seeds 1--10) used the package defaults: one initialization try, 2,000 iterations with the second half kept, and no tempered transitions (the defaults set the number of tempering rungs to one). A run takes about 9 seconds. Weakening the prior on the feature means (inverse-gamma$(2,20)$) did not change the default-start result, so the collapse is not driven by that prior.

\emph{Long runs} (seeds 1--5 of the same datasets) used 50,000 iterations with the second half kept, run in segments of 5,000 iterations with exact restarts from the saved state. The default start used the multiple-start algorithm with ten tries of 500 iterations, followed by \texttt{BMVMMM\_Theta\_est} with five tries. Tempering is available only in \texttt{BMVMMM\_warm\_start}, through a geometric ladder of $N_t$ rungs from inverse temperature 1 to $\beta_{N_t}$, applied once every given number of iterations. In test runs at $c=6$, the log acceptance probabilities of the first tempered transitions were between $-597$ and $-47$ for 50 rungs down to $\beta=0.3$, between $-368$ and $-89$ for 200 rungs down to $0.3$, and between $-262$ and $-1.8$ for 50 rungs down to $0.8$: with 2,000 latent memberships, any ladder hot enough to move between configurations is essentially never accepted, as \citet{marco2024functional} anticipate. We used 50 rungs down to $\beta=0.8$, one tempered transition every 500 iterations (100 per chain); acceptance rates were 0\% to 13\%. Split-$\hat R$ between the two starts was 12 to 26 for $s$ and 1.7 to 5.2 for $\log\alpha_3$. The two starts' 95\% intervals for $\alpha_3$ were disjoint on five of five datasets at $c=3$ and two of five at $c=6$, because the truth-started chains drift toward smaller $\alpha_3$. A chain took 150 to 200 seconds.

\paragraph{Reference posterior.}
The reference model is the same marginal model with memberships integrated out by quadrature on a logit grid with atoms at the endpoints. For speed, the sampler evaluates the likelihood on a 400-bin histogram of the data; every log-likelihood we report is recomputed exactly, without binning, with 128-node Gauss--Jacobi quadrature. Its seven global parameters are the two means, the effective covariance through a Cholesky factor, $\alpha_3$ and $\pi$. Both prior settings use $C_{\rm eff}=LL^\top+10^{-6}I$ with $L$ lower triangular, half-normal priors with scale 3 on the diagonal entries of $L$ (sampled on the log scale), $\mathcal N(0,3^2)$ on the off-diagonal entry, $\pi\sim\mathrm{Uniform}(0.02,0.98)$ (sampled on the logit scale) and $\alpha_3$ restricted to $[10^{-3},10^4]$. The weak setting puts $\mathcal N(0,10^2)$ on the means and $\mathcal N(\log4,1.5^2)$ on $\log\alpha_3$; the matched setting replaces these with \texttt{BayesFMMM}'s $t_{20}(0,0.316)$ and $\mathrm{Exp}(10)$, and does not reproduce the package's shrinkage prior on the covariance factors or its prior on $\sigma^2$. Each of 100 datasets (50 per value of $c$: the 20 of the sampler runs and 80 new ones) was run under both settings from two starts, the truth and the configuration in which \texttt{BayesFMMM}'s default start settles, with emcee using 28 walkers for 2,000 steps and the second half kept (400 runs; all runs seeded).

Under the matched setting the two starts agree on $\alpha_3$ (split-$\hat R\le1.15$), the posterior probability of $\alpha_3<0.1$ is at least 0.49 in every run, and the 95\% upper limits for $\alpha_3$ are 0.13 to 0.50. In 48 of the 200 runs more than 1\% (at most 3.5\%) of the draws of $\alpha_3$ fall below $1.5\times10^{-3}$, near the support bound; this affects the lower 2.5\% quantile, not the medians or upper limits. The hard posterior has two modes at $c=3$: nearly identical types ($s\approx0.1$ to $0.5$; 28 truth-start and 48 collapsed-start runs) and closer but distinct types ($s\approx1.0$ to $1.2$; the other 24 runs). At $c=6$ the types are nearly identical ($s\approx0.09$ to $0.25$) from both starts. The fit of the hard configuration is measured by the exact log-likelihood of the best posterior draw with $\alpha_3<0.1$ relative to the truth: medians $+2.03$ and $-0.71$ nats at $c=3$ (truth and collapsed starts) and $+1.66$ and $+1.64$ at $c=6$. Part of a positive value is the usual advantage of seven free parameters over the true values. We do not summarize by the coordinate-wise posterior median, which lies off the curved ridge of the posterior and fits 4 to 26 nats worse than the truth. An unbinned likelihood gives the same result on five datasets per value of $c$: the posterior is hard on all ten, and over all saved draws the binned and exact log-likelihoods relative to the truth differ by at most 0.15 nats.

Under the weak setting no run goes hard; every interval for $\alpha_3$ covers the truth, as do 192 of the 200 intervals for $s$. At $c=3$ the 2,000-step runs from the two starts mix incompletely (split-$\hat R$ up to 1.44); on the ten original datasets, 10,000 steps with 56 walkers give split-$\hat R\le1.07$ for both $\log\alpha_3$ and $s$, and every interval covers the truth.

\paragraph{Larger samples.}
We repeated the comparison at $c=3$ for $n\in\{10^4,\,3\times10^4,\,6\times10^4,\,1.2\times10^5\}$ and at $c=6$ for $n\in\{7.7\times10^4,\,3\times10^5,\,10^6\}$, five new datasets per $n$, with 2,000 histogram bins for $n\ge3\times10^4$ (400 and 2,000 bins agree within Monte Carlo error at $n=6\times10^4$). At these sizes the chains started at the collapsed configuration stay near it under both priors (split-$\hat R$ up to 2.3 for $\log\alpha_3$), so we report the truth-started chains. Under the package's priors the median of $\alpha_3$ rises only slowly with $n$: at $c=3$ it is about 0.05 up to $n=3\times10^4$, about 0.1 at $6\times10^4$ and 0.25 to 0.36 at $1.2\times10^5$; at $c=6$ and $n=10^6$ it is 1.0 to 1.8. The half-separation stays between 1.2 and 1.5, and no interval covers the true concentration or separation. The data do turn against the hard configuration: at $c=3$ the best hard draw goes from $+1.3$ nats relative to the truth at $n=2{,}000$ to $-5.9$ nats at $n=1.2\times10^5$, and at $c=6$, $n=10^6$ the chains no longer visit the hard region. Under the weak prior no run goes hard, and the truth-started intervals cover the truth in 77 of 80 cases. A rough accounting explains the persistence. Relative to a typical hard configuration (means $\pm0.15$, $\alpha_3=0.05$), the package's priors place the truth about 61 nats lower at $c=3$ and 90 nats lower at $c=6$, about 36 of them from the prior on the means, while the data's evidence against the closest hard look-alike grows like $n\KL$, with $\KL=1.5\times10^{-4}$ at $c=3$ and $5.9\times10^{-6}$ at $c=6$. The two balance only near $n\approx4\times10^5$ and $1.5\times10^7$. The certified sizes concern testing two fixed laws; they do not say when data overcome a prior that sits tens of nats away from the truth. In applications with 33 or 50 dimensions per unit, \citet{marco2024flexible} found their fits insensitive to the prior on $\alpha_3$; we expect the ambiguity, and with it the prior's leverage, to be weaker there (Section~\ref{sec:limits}).

Figure~\ref{fig:posterior} shows the four runs for dataset 1 at each value of $c$.

\section{Fitting the Fixed-Concentration Model}\label{app:em}

Each fit in Appendix~\ref{app:sweep} maximizes the marginal likelihood of the full model, with $c$ fixed, in two stages. An EM algorithm generates candidate starting points, and a quasi-Newton method then maximizes the marginal likelihood from the best of them. Only the second stage determines the reported fits, so the EM needs to explore well, not to be exact. All fits are reproducible: refits in fresh processes return identical values.

\subsection{Procedure}\label{app:em-code}

\paragraph{Data and parameters.}
The data are standardized by their sample mean and standard deviation, and the estimates are mapped back at the end. The free parameters are $\mu_1$, $\mu_2$, $\sigma_1$, $\sigma_2$ and the prior mean membership $\pi$.

\paragraph{Stage 1: starting points.}
The EM represents memberships on $M=256$ nodes placed at the quantiles $(m-\tfrac12)/M$ of $\Beta(c\pi,c(1-\pi))$, each with prior weight $1/M$, so the posterior weight of node $m$ for observation $i$ is proportional to $\varphi\big(x_i;\,p_m\mu_1+(1-p_m)\mu_2,\ p_m^2\sigma_1^2+(1-p_m)^2\sigma_2^2\big)$. The M-step alternates between the exact mean update of Proposition~\ref{prop:mstep-mu} and iteratively reweighted least squares for the variances with weights $g_m/V_m^2$ (Propositions~\ref{prop:mstep-sigma} and~\ref{prop:mstep-fisher}), and sets $\pi$ to the average posterior membership; a variance step that would make a variance negative is rejected and the previous value kept. The update of $\pi$ is a moment rule rather than the exact EM step, and the iteration does not check that the likelihood increases; both simplifications are harmless for an initializer. The EM stops when the relative change in its log-likelihood falls below $10^{-5}$, or after 50 iterations. It is run from the $k$-means solution (ten initializations, seeded by the dataset's seed) and from five perturbations of it: the means shifted by $\N(0,0.5^2)$ in standardized units, the log standard deviations by $\N(0,0.4^2)$, and $\pi$ drawn from a Dirichlet distribution with mean equal to the cluster proportions and concentration 10. For $c\ge6$ we also ran 39 perturbations.

\paragraph{Stage 2: maximization.}
From the three highest-likelihood EM endpoints, or for $c\ge6$ the ten best endpoints of the forty-start run whose separations differ by more than $10^{-2}$, we maximize the marginal log-likelihood directly with L-BFGS-B over $(\mu_1,\mu_2,\log\sigma_1,\log\sigma_2,\operatorname{logit}\pi)$, with $c$ fixed. The likelihood is computed with 128-node Gauss--Jacobi quadrature in the membership, with nodes and weights recomputed for the current $\Beta(c\pi,c(1-\pi))$ at every evaluation; for a fixed number of nodes the result is smooth in the parameters, and gradients are taken by finite differences. It agrees with 256 nodes to $1.2\times10^{-11}$ in relative terms at every fit. (The EM's quantile-midpoint rule is accurate only to about $10^{-5}$ at large $c$, enough for starting points but not for locating an optimum on a flat ridge.) The bounds, in standardized units, are $|\mu_k|\le20$, $\sigma_k\in[10^{-3},10]$ and $|\operatorname{logit}\pi|\le12$. A run that stops without success restarts from its own endpoint, at most twice, and the best maximized point is kept. This stage matters: at $n=10^4$ and $2\le c\le8$ it gained a median of 1.0 to 2.4 nats over the best EM endpoint, and up to 4.45. Nine runs, all at $c\le0.5$, ended on a line-search stop; in each, the three runs agreed to $10^{-3}$ in separation and within 0.01 nats. The largest gradient at a kept point corresponds to a displacement of order $10^{-5}$ in $\hat s$.

\paragraph{The likelihood has no global maximum.}
As for Gaussian mixtures with unequal variances, the likelihood is unbounded. Suppose the Beta exponent at a type's pure end is at most one; for type 1 this is $c(1-\pi)\le1$, which a free $\pi$ allows at any $c$. Put that type's mean at an observation and let its standard deviation tend to zero. The density at that observation then diverges, because near the pure end it behaves like $\int_0 y^{b-2}\,dy$ with $b=c(1-\pi)\le1$, while the densities at the other observations stay bounded below. The reported fits are therefore the highest interior local maxima reached from the starts, as is standard for mixtures. The point-mass fits of Appendix~\ref{app:sweep} are different: there the exponent exceeds one, the likelihood stays bounded, and it is flat as the standard deviation shrinks.

\subsection{The M-step of the initializer}\label{app:mstep}
For completeness we record the properties of the two M-step blocks.
Fix the nodes $p_m$ and the posterior weights $\gamma_{im}$, write $\sigma^2=(\sigma_1^2,\sigma_2^2)^\top$, $u^{(m)}=(p_m,1-p_m)^\top$, $q^{(m)}=(p_m^2,(1-p_m)^2)^\top$, $V_m=q^{(m)\top}\sigma^2$ and $e_i^{(m)}(\mu)=(x_i-u^{(m)\top}\mu)^2$. Up to constants, the expected complete-data log-likelihood in $(\mu,\sigma^2)$ is
\[
Q(\mu,\sigma^2)=-\tfrac12\sum_{i,m}\gamma_{im}\Big[\log V_m+\frac{e_i^{(m)}(\mu)}{V_m}\Big].
\]
It depends on the data only through the node sums $g_m=\sum_i\gamma_{im}$, $h_m=\sum_i\gamma_{im}x_i$ and $\rho_m=\sum_i\gamma_{im}x_i^2$, since $d_m(\mu):=\sum_i\gamma_{im}e_i^{(m)}(\mu)=\rho_m-2(u^{(m)\top}\mu)h_m+g_m(u^{(m)\top}\mu)^2$. Write $\bar e_m=d_m/g_m$ for the posterior-weighted mean squared residual at node $m$.

\begin{proposition}[Mean block]\label{prop:mstep-mu}
Fix $\sigma^2$ with positive entries. Then $\mu\mapsto Q$ is a concave quadratic with Hessian $-A$, where $A=\sum_m(g_m/V_m)\,u^{(m)}u^{(m)\top}$, and its stationary points solve $A\mu=b$ with $b=\sum_m(h_m/V_m)\,u^{(m)}$. Since $\gamma_{im}>0$ for all $i,m$, $A$ is positive definite as soon as two nodes are distinct, and $\mu=A^{-1}b$ is the unique maximizer of $Q$ in $\mu$.
\end{proposition}
\begin{proof}
$\partial_\mu Q=\sum_{i,m}\gamma_{im}u^{(m)}(x_i-u^{(m)\top}\mu)/V_m=b-A\mu$, and $z^\top Az=\sum_m(g_m/V_m)(z^\top u^{(m)})^2$ vanishes only if $z$ is orthogonal to every $u^{(m)}$.
\end{proof}

\begin{proposition}[Variance block: the weights are forced]\label{prop:mstep-sigma}
Fix $\mu$ and let $\sigma^2$ have positive entries. Then
\[
\frac{\partial Q}{\partial\sigma_k^2}=\tfrac12\sum_m q_k^{(m)}\,\frac{g_m(\bar e_m-V_m)}{V_m^2},
\]
so $\sigma^2$ is a stationary point of $Q$ if and only if it is a stationary point, at frozen weights, of the weighted least-squares criterion $\sum_m\omega_m(\bar e_m-q^{(m)\top}\sigma^2)^2$ with $\omega_m=g_m/V_m^2$. No other weight function has this property: if $\omega_m=g_m\psi(V_m)$ for a positive $\psi$, the normal equations agree with the stationarity condition for all data and all node configurations only if $\psi(V)\propto V^{-2}$.
\end{proposition}
\begin{proof}
$Q$ sees $\sigma^2$ only through $V_m$, and $\partial_V[-\tfrac12(g\log V+d/V)]=\tfrac12(d-gV)/V^2$; the chain rule and $d_m=g_m\bar e_m$ give the score. The normal equations of the weighted criterion are $\sum_m\omega_mq_k^{(m)}(\bar e_m-V_m)=0$, which coincide with the score equations exactly when $\omega_m\propto g_m/V_m^2$ with a constant common to all nodes; since the nodes realize a continuum of values of $V$, this forces $\psi(V)\propto V^{-2}$.
\end{proof}

The block is an exact Gamma regression: given $u$, $e(u;\mu)\sim V_\sigma(u)\chi^2_1$, whose log-density in its mean $V$ is $-\tfrac12\log V-e/(2V)$ up to terms free of $V$, and summing against $\gamma_{im}$ reproduces $Q$. The link is the identity, which is not canonical for the Gamma family, so $Q$ is not globally concave in $\sigma^2$.

\begin{proposition}[Where the block is concave]\label{prop:mstep-concave}
The Hessian of $Q$ in $\sigma^2$ is
\[
H=\tfrac12\sum_mq^{(m)}q^{(m)\top}g_m(V_m-2\bar e_m)/V_m^3 .
\]
Hence $Q$ is concave on the convex set $\{\sigma^2:V_m<2\bar e_m\text{ for every }m\}$, strictly so when the $q^{(m)}$ span $\R^2$, and then a stationary point in that set is its unique maximizer there.
\end{proposition}
\begin{proof}
Differentiate the score once more; on the set every factor $V_m-2\bar e_m$ is negative.
\end{proof}

\begin{proposition}[The variance update is Fisher scoring]\label{prop:mstep-fisher}
Let
\[
T(\sigma^2)=\Big[\sum_m\omega_mq^{(m)}q^{(m)\top}\Big]^{-1}\sum_m\omega_mq^{(m)}\bar e_m ,
\]
with $\omega_m=g_m/V_m^2$ evaluated at $\sigma^2$, let $\mathcal I=\tfrac12\sum_m\omega_mq^{(m)}q^{(m)\top}$ and let $U=\nabla_{\sigma^2}Q$. Then $T(\sigma^2)=\sigma^2+\mathcal I^{-1}U$. Its fixed points are the stationary points of $Q$ in $\sigma^2$, $T(\sigma^2)-\sigma^2$ is an ascent direction whenever $U\ne0$, and near a fixed point the iteration converges linearly at rate equal to the spectral radius of $\mathcal I^{-1}R$, with $R=\sum_mq^{(m)}q^{(m)\top}g_m(V_m-\bar e_m)/V_m^3$.
\end{proposition}
\begin{proof}
$\mathcal I\sigma^2+U=\tfrac12\sum_m\omega_mq^{(m)}V_m+\tfrac12\sum_m\omega_mq^{(m)}(\bar e_m-V_m)=\tfrac12\sum_m\omega_mq^{(m)}\bar e_m$; multiply by $\mathcal I^{-1}$. $\mathcal I$ is positive definite when the $q^{(m)}$ span $\R^2$. Since $H=-\mathcal I+R$, differentiating $T$ at a fixed point, where $U=0$, gives $DT=\mathcal I^{-1}R$.
\end{proof}

The rate is driven by the node-wise mismatch $\bar e_m-V_m$, which has conditional mean zero at a correctly specified configuration, so the iteration is fast where the model fits. A unit scoring step is not guaranteed to keep the variances positive; the implementation rejects such steps (Appendix~\ref{app:em-code}).

%% file: sections/appendix_proofs.tex
\section{Proofs for Section~\ref{sec:ident}}\label{app:ident}

\subsection{The symmetric model}
The symmetric model needs only a one-term expansion, which also covers the case $s=0$ excluded from Proposition~\ref{prop:ident}.

\begin{lemma}[Tail of the moment generating function]\label{lem:mgf-sym}
Let $s\ge0$, $r>0$, $t>0$, write $M_\theta(z)=\mathbb E\,e^{zX}$ and $A_\theta(z)=sz+t^2z^2/2$. As $z\to+\infty$,
\[
\log M_\theta(z)=A_\theta(z)-r\log A_\theta(z)+\log\frac{\Gamma(2r)}{2^r\,\Gamma(r)}+\mathbb 1\{s=0\}\log2+o(1).
\]
\end{lemma}
\begin{proof}
$M_\theta$ is finite on $\mathbb R$ because $|W|\le1$. The density of $W$ is $p_r(w)=C_r(1-w^2)^{r-1}$ on $[-1,1]$ with $C_r=\Gamma(2r)/\big(2^{2r-1}\Gamma(r)^2\big)$, and conditioning on $W$ gives $M_\theta(z)=\mathbb E\exp\{\varphi_z(W)\}$ with $\varphi_z(w)=szw+\tfrac14t^2z^2(1+w^2)$. Substituting $w=1-y$ gives $\varphi_z(1-y)=A-szy-\tfrac14t^2z^2\,y(2-y)$, where $A=A_\theta(z)$, so
\[
M_\theta(z)=C_r\,e^{A}\,J(z),\qquad J(z)=\int_0^2 e^{-szy-\frac14t^2z^2y(2-y)}\,\big(y(2-y)\big)^{r-1}\,dy .
\]
Fix $\delta\in(0,1)$ and split $J=J_L+J_M+J_R$ over $[0,\delta]$, $[\delta,2-\delta]$, $[2-\delta,2]$.
\emph{Middle.} There $y(2-y)\ge\delta$, so $J_M\le e^{-t^2z^2\delta/4}\int_0^2(y(2-y))^{r-1}dy=O(e^{-t^2z^2\delta/4})$.
\emph{Left.} On $[0,\delta]$, $2y(1-\delta/2)\le y(2-y)\le2y$, so the exponent lies between $-Ay$ and $-(1-\delta/2)Ay$; also $(y(2-y))^{r-1}=(2y)^{r-1}(1-y/2)^{r-1}$ with $(1-y/2)^{r-1}\in[m_\delta,M_\delta]$, where $m_\delta,M_\delta\to1$ as $\delta\to0$. Since $\int_0^\delta e^{-By}y^{r-1}dy=\Gamma(r)B^{-r}(1+o(1))$ as $B\to\infty$,
\[
m_\delta\,2^{r-1}\Gamma(r)\le\liminf_{z\to\infty}A^rJ_L\le\limsup_{z\to\infty}A^rJ_L\le M_\delta\,2^{r-1}\Gamma(r)(1-\delta/2)^{-r}.
\]
\emph{Right.} With $y=2-u$, the exponent is $-2sz+szu-\tfrac14t^2z^2u(2-u)\le-2sz-B'u$ on $[0,\delta]$, where $B'=(1-\delta/2)t^2z^2/2-sz\to\infty$. Hence $A^rJ_R=O(e^{-2sz})\to0$ if $s>0$. If $s=0$, the integrand is invariant under $y\mapsto2-y$, so $J_R=J_L$.
Letting $z\to\infty$ and then $\delta\to0$ gives $A^rJ(z)\to2^{r-1}\Gamma(r)\,(1+\mathbb 1\{s=0\})$. Finally, $C_r\,2^{r-1}\Gamma(r)=\Gamma(2r)/(2^r\Gamma(r))$.
\end{proof}

\begin{proposition}[Identifiability, symmetric model]\label{prop:ident-sym}
For $\theta=(s,r,t)$ and $\theta'=(s',r',t')$ in $\mathbb R\times(0,\infty)\times(0,\infty)$, $f_\theta=f_{\theta'}$ if and only if $(|s|,r,t)=(|s'|,r',t')$.
\end{proposition}
\begin{proof}
Because $W\overset{d}{=}-W$, we have $f_{(s,r,t)}=f_{(-s,r,t)}$, so it suffices to take $s,s'\ge0$. If $f_\theta=f_{\theta'}$, then $M_\theta=M_{\theta'}$ on $\mathbb R$. By Lemma~\ref{lem:mgf-sym}, each of the following limits depends on $M_\theta$ alone:
$t^2=\lim_{z\to\infty}2\log M_\theta(z)/z^2$;
$s=\lim_{z\to\infty}\big(\log M_\theta(z)-t^2z^2/2\big)/z$, since $r\log A_\theta(z)=O(\log z)$;
$r=-\tfrac12\lim_{z\to\infty}\big(\log M_\theta(z)-t^2z^2/2-sz\big)/\log z$, since $\log A_\theta(z)=2\log z+\log(t^2/2)+o(1)$.
Hence $(s,r,t)=(s',r',t')$.
\end{proof}

\subsection{The two-type model}
\begin{proof}[Proof of Lemma~\ref{lem:mgf2}]
Given $U=u$, $X\sim\mathcal N\big(\mu_2+u\Delta,\ \sigma_1^2u^2+\sigma_2^2(1-u)^2\big)$, so, with $p(u)=u^{a-1}(1-u)^{b-1}/B(a,b)$ the density of $U$, $M(z)=\mathbb E\,e^{\varphi_z(U)}$ with $\varphi_z(u)=z(\mu_2+u\Delta)+\tfrac12z^2\big(\sigma_1^2u^2+\sigma_2^2(1-u)^2\big)$. Put $T(z)=\mu_1z+\tfrac12\sigma_1^2z^2=\varphi_z(1)$, $B_z=\sigma_1^2z^2+\Delta z$, $q=\tfrac12(\sigma_1^2+\sigma_2^2)\le\sigma_1^2$, and $\delta=\tfrac12$.

\emph{Away from $u=1$.} For $u\in[0,1-\delta]$, $\varphi_z(u)-T(z)=-(1-u)\big[z\Delta+\tfrac12z^2\big(\sigma_1^2(1+u)-\sigma_2^2(1-u)\big)\big]$. If $\sigma_1>\sigma_2$, the bracket is at least $z\Delta+\tfrac12z^2(\sigma_1^2-\sigma_2^2)\ge\tfrac14(\sigma_1^2-\sigma_2^2)z^2$ for large $z$. If $\sigma_1=\sigma_2$, it is at least $z\Delta$ with $\Delta>0$. Either way $\int_0^{1-\delta}e^{\varphi_z(u)-T(z)}\,p(u)\,du\le e^{-cz}$ for some $c>0$ and large $z$, which is $O(z^{-k})$ relative to $B_z^{-b}$ for every $k$.

\emph{Near $u=1$.} With $y=1-u$, $\varphi_z(1-y)-T(z)=-B_zy+qz^2y^2$ and $p(1-y)=y^{b-1}(1-y)^{a-1}/B(a,b)$. Substituting $y=\tau/B_z$,
\[
\int_0^{\delta}e^{-B_zy+qz^2y^2}y^{b-1}(1-y)^{a-1}dy
=B_z^{-b}\int_0^{\delta B_z}e^{-\tau}\tau^{b-1}g_z(\tau)\,d\tau,\qquad
g_z(\tau)=e^{\kappa_z\tau^2}\Big(1-\frac{\tau}{B_z}\Big)^{a-1},
\]
with $\kappa_z=qz^2/B_z^2=O(z^{-2})$ and $1/B_z=O(z^{-2})$.
On $\tau\le\sqrt z$, Taylor expansion gives $g_z(\tau)=1+\kappa_z\tau^2-(a-1)\tau/B_z+R_z(\tau)$ with $|R_z(\tau)|\le Cz^{-4}(1+\tau^4)$.
On $\sqrt z<\tau\le\delta B_z$, $\kappa_z\tau\le q\delta z^2/B_z\to q\delta/\sigma_1^2\le\tfrac12$, so $e^{-\tau+\kappa_z\tau^2}\le e^{-\tau/2}$ for large $z$, and $(1-\tau/B_z)^{a-1}\le\max(1,2^{1-a})$; this range contributes $O(e^{-\sqrt z/4})$, as do the tails beyond $\sqrt z$ of the integrals of the Taylor terms. Therefore
\[
\int_0^{\delta B_z}e^{-\tau}\tau^{b-1}g_z(\tau)\,d\tau=\Gamma(b)\Big[1+\kappa_z\,b(b+1)-\frac{(a-1)b}{B_z}\Big]+O(z^{-4})
=\Gamma(b)\Big[1+\frac{K}{z^2}+O(z^{-3})\Big],
\]
with $K=\frac{b(b+1)(\sigma_1^2+\sigma_2^2)}{2\sigma_1^4}-\frac{(a-1)b}{\sigma_1^2}$, using $\kappa_z=\frac{q}{\sigma_1^4z^2}+O(z^{-3})$ and $B_z^{-1}=\frac{1}{\sigma_1^2z^2}+O(z^{-3})$.

\emph{Assembly.} $\log M(z)=T(z)-b\log B_z+\log\frac{\Gamma(b)}{B(a,b)}+\frac{K}{z^2}+O(z^{-3})$, and $-b\log B_z=-2b\log z-b\log\sigma_1^2-\frac{b\Delta}{\sigma_1^2z}+\frac{b\Delta^2}{2\sigma_1^4z^2}+O(z^{-3})$.
\end{proof}

\paragraph{Numerical check.} The expansion was checked with mpmath at 50 digits for three parameter sets, including $\mu_1=\mu_2$ with $\sigma_1\ne\sigma_2$. After subtracting the expansion through orders $0$, $1$ and $2$, the residuals shrink by factors of $2$, $4$ and $8$ per doubling of $z$ (by $16$ at order $2$ when $\Delta=0$, where the odd terms vanish).

\begin{proof}[Proof of Proposition~\ref{prop:ident}]
The statement compares parameter vectors within the class $(\mu_1,\sigma_1)\ne(\mu_2,\sigma_2)$. Whether a law in this class can also arise from the excluded case, where both endpoints of $[0,1]$ contribute to the tail of $M$, is not settled by the expansion, and we do not claim it; within the symmetric family, Proposition~\ref{prop:ident-sym} covers that case ($s=0$). The label swap preserves the law because $1-U\sim\mathrm{Beta}(b,a)$, and exactly one of the two labelings satisfies $(\sigma_1,\mu_1)\succ(\sigma_2,\mu_2)$. Let $\vartheta,\vartheta'$ both satisfy this ordering and induce the same law. Then $M=M'$ on $\mathbb R$, and since $z^2,z,\log z,1,z^{-1},z^{-2}$ form an asymptotic scale, the coefficients in Lemma~\ref{lem:mgf2} agree:
\begin{enumerate}
\item the $z^2$ coefficient gives $\sigma_1=\sigma_1'$;
\item the $z$ coefficient gives $\mu_1=\mu_1'$;
\item the $\log z$ coefficient gives $b=b'$;
\item $C_0=C_0'$ gives $B(a,b)=B(a',b)$, hence $a=a'$, because $\partial_a\log B(a,b)=\psi(a)-\psi(a+b)<0$;
\item $C_1=C_1'$ gives $\Delta=\Delta'$, hence $\mu_2=\mu_2'$;
\item $C_2=C_2'$ gives $\sigma_2^2=\sigma_2'^2$, because $\sigma_2^2$ enters $C_2$ linearly with coefficient $b(b+1)/(2\sigma_1^4)>0$ and every other term already agrees.
\end{enumerate}
\end{proof}

\section{Proofs for Sections~\ref{sec:hard} and~\ref{sec:gauss}}\label{app:hard}

\begin{proof}[Proof of Lemma~\ref{lem:cumulants}]
$W$ is symmetric and $W^2\sim\mathrm{Beta}(\tfrac12,r)$, so $m_2:=\mathbb E W^2=1/(2r+1)$ and $m_4:=\mathbb E W^4=3/\big((2r+1)(2r+3)\big)$. Conditioning on $W$, with $\mu=sW$ and $\sigma^2=t^2(1+W^2)/2$,
$\mathbb E X^2=\mathbb E[\mu^2+\sigma^2]$ and $\mathbb E X^4=\mathbb E[\mu^4+6\mu^2\sigma^2+3\sigma^4]$. Substituting gives
$\kappa_4=\mathbb E X^4-3v^2 = s^4(m_4-3m_2^2)+3s^2t^2(m_4-m_2^2)+\tfrac34 t^4(m_4-m_2^2)$, and
$m_4-m_2^2 = 4r/\big((2r+1)^2(2r+3)\big)$, $m_4-3m_2^2=-6/\big((2r+1)^2(2r+3)\big)$. Odd cumulants vanish by the symmetry $X\overset{d}{=}-X$.
\end{proof}

\begin{proof}[Proof of Proposition~\ref{prop:range}]
Lower bound: $\mathbb E X^4\ge(\mathbb E X^2)^2$ with equality only if $|X|$ is a.s.\ constant, which is impossible since $f_\theta$ is a density; hence $\gamma_2=\mathbb E X^4/v^2-3>-2$.
Upper bound: drop the $-6s^4$ term and use $(s^2+(r+1)t^2)^2\ge\max\{4(r+1)s^2t^2,\,(r+1)^2t^4\}$ on the two remaining terms separately:
$\gamma_2 \le \frac{3r}{(2r+3)(r+1)}+\frac{3r}{(2r+3)(r+1)^2} < \tfrac32+\tfrac32 = 3$.
(Numerically, $\sup_\theta\gamma_2\approx0.24$.)
Sign: the numerator of $\kappa_4$ is negative iff $r(12s^2t^2+3t^4)<6s^4$, i.e.\ $r<2d^4/(4d^2+1)$.
\end{proof}

\begin{proof}[Proof of Proposition~\ref{prop:lookalike}]
Well-definedness is Proposition~\ref{prop:range}: $t_G^2>0\iff\gamma_2>-2$ and $v(1-\sqrt{\gamma_2/3})>0\iff\gamma_2<3$.
Case $\gamma_2<0$: write $Y=s_G B+t_G Z$ with $B$ uniform on $\{\pm1\}$ independent of $Z$. Cumulants add over independent summands, $\kappa_k(t_GZ)=0$ for $k\ge3$, and $B$ has $\kappa_2=1$, $\kappa_4=\mathbb E B^4-3=-2$ and vanishing odd cumulants. Hence $\kappa_2(Y)=s_G^2+t_G^2=v$ and $\kappa_4(Y)=-2s_G^4=\gamma_2v^2=\kappa_4(\theta)$; odd cumulants vanish on both sides.
Case $\gamma_2>0$: $Y=\sigma Z$ with $\sigma^2$ uniform on $\{v(1\pm\sqrt{\gamma_2/3})\}$ has $\kappa_2=\mathbb E\sigma^2=v$ and $\kappa_4=3\operatorname{Var}(\sigma^2)=\gamma_2v^2$.
The limit statement follows from $W\to\pm1$ in distribution as $r\to0$ and $(1+W^2)/2\to1$.
\end{proof}

\begin{proof}[Proof of Theorem~\ref{thm:lecam}]
The minimal error sum over tests equals $1-\TV(f_\theta^{\otimes n},G_\theta^{\otimes n})$. By Pinsker's inequality and additivity of $\KL$ over products, $\TV(f_\theta^{\otimes n},G_\theta^{\otimes n})\le\sqrt{n\,\KL(f_\theta\|G_\theta)/2}=\tfrac12\sqrt{n/\nhard(\theta)}$. $\nhard$ is positive because $\KL(f_\theta\|G_\theta)<\infty$: $\log(f_\theta/G_\theta)$ grows at most quadratically in $x$, because both are finite or continuous Gaussian location-scale mixtures with variances bounded away from zero, and $f_\theta$ has Gaussian tails. $\nhard$ is finite because $\KL(f_\theta\|G_\theta)>0$, that is, $f_\theta\ne G_\theta$: by Lemma~\ref{lem:mgf-sym}, $\log M_\theta(z)$ contains the term $-r\log A_\theta(z)\sim-2r\log z$ as $z\to\infty$, whereas the log moment generating function of every candidate $G_\theta$ is a quadratic in $z$ plus a bounded term ($t_G^2z^2/2+s_G z-\log2+o(1)$, $vz^2/2$, or $v_+z^2/2-\log2+o(1)$). Dependence on $(d,c)$ only follows from the scaling identity in Section~\ref{sec:model}, since $G_\theta$ is built from $(v,\gamma_2)$ in a scale-equivariant way.
\end{proof}

\begin{proposition}[Tightness for a fixed pair]\label{prop:tight}
Let $H^2=H^2(f_\theta,G_\theta)=\tfrac12\int(\sqrt{f_\theta}-\sqrt{G_\theta})^2$. Every test $\phi$ based on $n$ observations satisfies
\[
f_\theta^{\otimes n}(\phi=1)+G_\theta^{\otimes n}(\phi=0)\ge\tfrac12\big(1-H^2\big)^{2n},
\]
and the likelihood-ratio test $\phi^\star=\mathbb 1\{\prod_iG_\theta(X_i)>\prod_if_\theta(X_i)\}$ satisfies
\[
f_\theta^{\otimes n}(\phi^\star=1)+G_\theta^{\otimes n}(\phi^\star=0)\le\big(1-H^2\big)^n .
\]
Consequently, the smallest $n$ at which some test has summed error at most $\epsilon<\tfrac12$ lies between $\log(1/(2\epsilon))/\{-2\log(1-H^2)\}$ and $\log(1/\epsilon)/\{-\log(1-H^2)\}$, which to first order in $H^2$ are $\log(1/(2\epsilon))/(2H^2)$ and $\log(1/\epsilon)/H^2$.
\end{proposition}
\begin{proof}
The minimal error sum over tests is $1-\TV(f_\theta^{\otimes n},G_\theta^{\otimes n})=\int\min(f_\theta^{\otimes n},G_\theta^{\otimes n})$. Le Cam's inequality \citep[Lemma~2.3]{tsybakov2009} bounds it below by $\tfrac12\big(\int\sqrt{f_\theta^{\otimes n}G_\theta^{\otimes n}}\big)^2=\tfrac12(1-H^2)^{2n}$. The error sum of $\phi^\star$ equals $\int\min(f_\theta^{\otimes n},G_\theta^{\otimes n})\le\int\sqrt{f_\theta^{\otimes n}G_\theta^{\otimes n}}=(1-H^2)^n$. Solving each bound for $n$ at error $\epsilon$ gives the two ends.
\end{proof}

At $\epsilon=0.05$ the two ends differ by the factor $\log 20/(\tfrac12\log 10)\approx2.6$, whatever $H^2$. Numerically, for the closed-form look-alike at $d=3$, $H^2/\KL\in[0.250,0.252]$ for $c\in\{1,2,3,4,6,8\}$, so the two ends are about $9.1\,\nhard$ and $24\,\nhard$ (Table~\ref{tab:tight}).

\begin{lemma}[Continuity at hard membership]\label{lem:continuity}
For $\gamma_2(\theta)<0$ and $\theta_\epsilon=(s_G,\epsilon/2,t_G)$, $\KL(f_\theta\|f_{\theta_\epsilon})\to\KL(f_\theta\|G_\theta)$ as $\epsilon\downarrow0$.
\end{lemma}
\begin{proof}
Under $\theta_\epsilon$, given $U=u$ the density of $X$ at $x$ is $g_x(u)=\varphi\big(x;\,s_G(2u-1),\,t_G^2(u^2+(1-u)^2)\big)$, with mean in $[-s_G,s_G]$ and variance in $[t_G^2/2,\,t_G^2]$. For fixed $x$, $u\mapsto g_x(u)$ is bounded and continuous on $[0,1]$, and $\mathrm{Beta}(\epsilon/2,\epsilon/2)\Rightarrow\tfrac12\delta_0+\tfrac12\delta_1$, so $f_{\theta_\epsilon}(x)\to\tfrac12g_x(0)+\tfrac12g_x(1)=G_\theta(x)$ pointwise. The variance bounds give, uniformly in $\epsilon$,
$(2\pi t_G^2)^{-1/2}\exp\{-(|x|+s_G)^2/t_G^2\}\le f_{\theta_\epsilon}(x)\le(\pi t_G^2)^{-1/2}$,
so $|\log f_{\theta_\epsilon}(x)|\le C+(|x|+s_G)^2/t_G^2$, which is $f_\theta$-integrable since $f_\theta$ has a finite second moment. Dominated convergence gives $\int f_\theta\log f_{\theta_\epsilon}\to\int f_\theta\log G_\theta$; $\int f_\theta\log f_\theta$ is finite.
\end{proof}

\begin{proof}[Proof of Corollary~\ref{cor:estimation}]
Le Cam's two-point reduction: for $\theta_0,\theta_1$ with $|c(\theta_0)-c(\theta_1)|\ge2\delta$, any estimator has $\max_i\mathbb P_{\theta_i}(|\hat c-c(\theta_i)|\ge\delta)\ge\tfrac12\{1-\TV(f_{\theta_0}^{\otimes n},f_{\theta_1}^{\otimes n})\}$. Take $\theta_0=\theta$, $\theta_1=(s_G,\epsilon/2,t_G)$ and $\delta=(c-\epsilon)/2$. By Pinsker, $\TV\le\sqrt{n\KL(f_\theta\|f_{\theta_1})/2}$, and by Lemma~\ref{lem:continuity} this tends to $\sqrt{n\KL(f_\theta\|G_\theta)/2}=\tfrac12\sqrt{n/\nhard}\le\tfrac12$ as $\epsilon\to0$.
\end{proof}

\begin{proof}[Proof of Proposition~\ref{prop:gauss-end}]
Write $s_k=s_{r_k}$. Since $t^2<2v$, $s_k^2=r_k(2v-t^2)+(v-t^2)>0$ for large $k$. Let $W_k=2U_k-1$ with $U_k\sim\mathrm{Beta}(r_k,r_k)$. By the central limit theorem for symmetric Beta laws, $\sqrt{2r_k+1}\,W_k\Rightarrow\mathcal N(0,1)$, and $s_k^2/(2r_k+1)\to v-t^2/2$, so $s_kW_k\Rightarrow\mathcal N(0,v-t^2/2)$ while $W_k^2\to0$ in probability. By Slutsky's lemma, $(s_kW_k,W_k^2)\Rightarrow(N,0)$ jointly, with $N\sim\mathcal N(0,v-t^2/2)$. For fixed $x$, the map $(m,w)\mapsto\varphi\big(x;m,t^2(1+w)/2\big)$ is bounded and continuous on $\mathbb R\times[0,1]$, so by the definition of weak convergence (the portmanteau theorem)
$f_{\theta_k}(x)=\mathbb E\,\varphi\big(x;s_kW_k,t^2(1+W_k^2)/2\big)\to\mathbb E\,\varphi(x;N,t^2/2)=\varphi(x;0,v)$. Scheff\'e's lemma upgrades pointwise convergence of densities to total variation. Finally,
$\TV(f_\theta^{\otimes n},f_{\theta_k}^{\otimes n})\le\TV(f_\theta^{\otimes n},\mathcal N(0,v)^{\otimes n})+n\,\TV(\mathcal N(0,v),f_{\theta_k})$; the first term is bounded by Pinsker as in Theorem~\ref{thm:lecam}, and the second tends to $0$ for fixed $n$.
\end{proof}

\begin{proof}[Proof of Proposition~\ref{prop:scaling}]
From Lemma~\ref{lem:cumulants}, $\gamma_2=\dfrac{-6s^4+12rs^2t^2+3rt^4}{(2r+3)\big(s^2+(r+1)t^2\big)^2}$.
(a) Substitute $s^2=r\,v(2-\rho)+v(1-\rho)$: the numerator is $6r^2v^2(2-\rho)(3\rho-2)+O(r)$ and the denominator is $8r^3v^2+O(r^2)$.
(b) The numerator is $r(12s^2t^2+3t^4)+O(1)$ and the denominator is $2r^3t^4+O(r^2)$.
\end{proof}

%% file: main.bbl
\begin{thebibliography}{30}
\providecommand{\natexlab}[1]{#1}
\providecommand{\url}[1]{\texttt{#1}}
\expandafter\ifx\csname urlstyle\endcsname\relax
  \providecommand{\doi}[1]{doi: #1}\else
  \providecommand{\doi}{doi: \begingroup \urlstyle{rm}\Url}\fi

\bibitem[Airoldi et~al.(2008)Airoldi, Blei, Fienberg, and Xing]{airoldi2008}
Edoardo~M. Airoldi, David~M. Blei, Stephen~E. Fienberg, and Eric~P. Xing.
\newblock Mixed membership stochastic blockmodels.
\newblock \emph{Journal of Machine Learning Research}, 9:\penalty0 1981--2014,
  2008.

\bibitem[Andrews and Cheng(2012)]{andrews2012}
Donald W.~K. Andrews and Xu~Cheng.
\newblock Estimation and inference with weak, semi-strong, and strong
  identification.
\newblock \emph{Econometrica}, 80\penalty0 (5):\penalty0 2153--2211, 2012.

\bibitem[Blei et~al.(2003)Blei, Ng, and Jordan]{blei2003}
David~M. Blei, Andrew~Y. Ng, and Michael~I. Jordan.
\newblock Latent {D}irichlet allocation.
\newblock \emph{Journal of Machine Learning Research}, 3:\penalty0 993--1022,
  2003.

\bibitem[Chen and Chen(2001)]{chen2001}
Hanfeng Chen and Jiahua Chen.
\newblock The likelihood ratio test for homogeneity in finite mixture models.
\newblock \emph{The Canadian Journal of Statistics}, 29\penalty0 (2):\penalty0
  201--215, 2001.

\bibitem[Eches et~al.(2010)Eches, Dobigeon, Mailhes, and Tourneret]{eches2010}
Olivier Eches, Nicolas Dobigeon, Corinne Mailhes, and Jean-Yves Tourneret.
\newblock Bayesian estimation of linear mixtures using the normal compositional
  model. application to hyperspectral imagery.
\newblock \emph{IEEE Transactions on Image Processing}, 19\penalty0
  (6):\penalty0 1403--1413, 2010.

\bibitem[Eismann and Stein(2007)]{eismann2007}
Michael~T. Eismann and David W.~J. Stein.
\newblock Stochastic mixture modeling.
\newblock In Chein-I Chang, editor, \emph{Hyperspectral Data Exploitation:
  Theory and Applications}, pages 107--148. Wiley, 2007.

\bibitem[Erosheva et~al.(2004)Erosheva, Fienberg, and Lafferty]{erosheva2004}
Elena Erosheva, Stephen Fienberg, and John Lafferty.
\newblock Mixed-membership models of scientific publications.
\newblock \emph{Proceedings of the National Academy of Sciences}, 101\penalty0
  (suppl 1):\penalty0 5220--5227, 2004.

\bibitem[Foreman-Mackey et~al.(2013)Foreman-Mackey, Hogg, Lang, and
  Goodman]{foreman2013}
Daniel Foreman-Mackey, David~W. Hogg, Dustin Lang, and Jonathan Goodman.
\newblock emcee: The {MCMC} hammer.
\newblock \emph{Publications of the Astronomical Society of the Pacific},
  125\penalty0 (925):\penalty0 306--312, 2013.

\bibitem[Ghosh and Sen(1985)]{ghosh1985}
Jayanta~K. Ghosh and Pranab~K. Sen.
\newblock On the asymptotic performance of the log likelihood ratio statistic
  for the mixture model and related results.
\newblock In Lucien~M. Le~Cam and Richard~A. Olshen, editors, \emph{Proceedings
  of the Berkeley Conference in Honor of Jerzy Neyman and Jack Kiefer},
  volume~2, pages 789--806. Wadsworth, 1985.

\bibitem[Gruhl and Erosheva(2014)]{gruhl2014}
Jonathan Gruhl and Elena~A. Erosheva.
\newblock A tale of two (types of) memberships: Comparing mixed and partial
  membership with a continuous data example.
\newblock In Edoardo~M. Airoldi, David~M. Blei, Elena~A. Erosheva, and
  Stephen~E. Fienberg, editors, \emph{Handbook of Mixed Membership Models and
  Their Applications}, pages 15--38. Chapman \& Hall/CRC, 2014.

\bibitem[Gu et~al.(2023)Gu, Erosheva, Xu, and Dunson]{gu2023}
Yuqi Gu, Elena~A. Erosheva, Gongjun Xu, and David~B. Dunson.
\newblock Dimension-grouped mixed membership models for multivariate
  categorical data.
\newblock \emph{Journal of Machine Learning Research}, 24\penalty0
  (88):\penalty0 1--49, 2023.

\bibitem[Halimi et~al.(2015)Halimi, Dobigeon, and Tourneret]{halimi2015}
Abderrahim Halimi, Nicolas Dobigeon, and Jean-Yves Tourneret.
\newblock Unsupervised unmixing of hyperspectral images accounting for
  endmember variability.
\newblock \emph{IEEE Transactions on Image Processing}, 24\penalty0
  (12):\penalty0 4904--4917, 2015.

\bibitem[Hardt and Price(2015)]{hardt2015}
Moritz Hardt and Eric Price.
\newblock Tight bounds for learning a mixture of two {G}aussians.
\newblock In \emph{Proceedings of the 47th Annual ACM Symposium on Theory of
  Computing (STOC)}, pages 753--760, 2015.

\bibitem[Heinrich and Kahn(2018)]{heinrich2018}
Philippe Heinrich and Jonas Kahn.
\newblock Strong identifiability and optimal minimax rates for finite mixture
  estimation.
\newblock \emph{The Annals of Statistics}, 46\penalty0 (6A):\penalty0
  2844--2870, 2018.

\bibitem[Heller et~al.(2008)Heller, Williamson, and Ghahramani]{heller2008}
Katherine~A. Heller, Sinead Williamson, and Zoubin Ghahramani.
\newblock Statistical models for partial membership.
\newblock In \emph{Proceedings of the 25th International Conference on Machine
  Learning}, pages 392--399, 2008.

\bibitem[Ho and Nguyen(2016)]{ho2016}
Nhat Ho and XuanLong Nguyen.
\newblock Convergence rates of parameter estimation for some weakly
  identifiable finite mixtures.
\newblock \emph{The Annals of Statistics}, 44\penalty0 (6):\penalty0
  2726--2755, 2016.

\bibitem[Le~Cam(1973)]{lecam1973}
Lucien Le~Cam.
\newblock Convergence of estimates under dimensionality restrictions.
\newblock \emph{The Annals of Statistics}, 1\penalty0 (1):\penalty0 38--53,
  1973.

\bibitem[Marco et~al.(2024{\natexlab{a}})Marco, {\c{S}}ent{\"u}rk, Jeste,
  DiStefano, Dickinson, and Telesca]{marco2024functional}
Nicholas Marco, Damla {\c{S}}ent{\"u}rk, Shafali Jeste, Charlotte DiStefano,
  Abigail Dickinson, and Donatello Telesca.
\newblock Functional mixed membership models.
\newblock \emph{Journal of Computational and Graphical Statistics}, 33\penalty0
  (4):\penalty0 1139--1149, 2024{\natexlab{a}}.
\newblock \doi{10.1080/10618600.2024.2304633}.

\bibitem[Marco et~al.(2024{\natexlab{b}})Marco, {\c{S}}ent{\"u}rk, Jeste,
  DiStefano, Dickinson, and Telesca]{marco2024flexible}
Nicholas Marco, Damla {\c{S}}ent{\"u}rk, Shafali Jeste, Charlotte~C. DiStefano,
  Abigail Dickinson, and Donatello Telesca.
\newblock Flexible regularized estimation in high-dimensional mixed membership
  models.
\newblock \emph{Computational Statistics \& Data Analysis}, 194:\penalty0
  107931, 2024{\natexlab{b}}.
\newblock \doi{10.1016/j.csda.2024.107931}.

\bibitem[Nguyen(2013)]{nguyen2013}
XuanLong Nguyen.
\newblock Convergence of latent mixing measures in finite and infinite mixture
  models.
\newblock \emph{The Annals of Statistics}, 41\penalty0 (1):\penalty0 370--400,
  2013.

\bibitem[Pritchard et~al.(2000)Pritchard, Stephens, and
  Donnelly]{pritchard2000}
Jonathan~K. Pritchard, Matthew Stephens, and Peter Donnelly.
\newblock Inference of population structure using multilocus genotype data.
\newblock \emph{Genetics}, 155\penalty0 (2):\penalty0 945--959, 2000.

\bibitem[Stein(2003)]{stein2003}
David Stein.
\newblock Application of the normal compositional model to the analysis of
  hyperspectral imagery.
\newblock In \emph{IEEE Workshop on Advances in Techniques for Analysis of
  Remotely Sensed Data}, pages 44--51, 2003.

\bibitem[Teicher(1963)]{teicher1963}
Henry Teicher.
\newblock Identifiability of finite mixtures.
\newblock \emph{The Annals of Mathematical Statistics}, 34\penalty0
  (4):\penalty0 1265--1269, 1963.

\bibitem[Tsybakov(2009)]{tsybakov2009}
Alexandre~B. Tsybakov.
\newblock \emph{Introduction to Nonparametric Estimation}.
\newblock Springer, New York, 2009.

\bibitem[Watanabe(2009)]{watanabe2009}
Sumio Watanabe.
\newblock \emph{Algebraic Geometry and Statistical Learning Theory}.
\newblock Cambridge University Press, 2009.

\bibitem[Woodbridge et~al.(2019)Woodbridge, Okun, Elidan, and
  Wiesel]{woodbridge2019}
Yonatan Woodbridge, Uri Okun, Gal Elidan, and Ami Wiesel.
\newblock Unmixing {$K$}-{G}aussians with application to hyperspectral imaging.
\newblock \emph{IEEE Transactions on Geoscience and Remote Sensing},
  57\penalty0 (9):\penalty0 7281--7293, 2019.

\bibitem[Wu and Yang(2020)]{wu2020}
Yihong Wu and Pengkun Yang.
\newblock Optimal estimation of {G}aussian mixtures via denoised method of
  moments.
\newblock \emph{The Annals of Statistics}, 48\penalty0 (4):\penalty0
  1981--2007, 2020.

\bibitem[Yakowitz and Spragins(1968)]{yakowitz1968}
Sidney~J. Yakowitz and John~D. Spragins.
\newblock On the identifiability of finite mixtures.
\newblock \emph{The Annals of Mathematical Statistics}, 39\penalty0
  (1):\penalty0 209--214, 1968.

\bibitem[Yurochkin et~al.(2019)Yurochkin, Guha, Sun, and Nguyen]{yurochkin2019}
Mikhail Yurochkin, Aritra Guha, Yuekai Sun, and XuanLong Nguyen.
\newblock {D}irichlet simplex nest and geometric inference.
\newblock In \emph{Proceedings of the 36th International Conference on Machine
  Learning}, volume~97 of \emph{PMLR}, pages 7262--7271, 2019.

\bibitem[Zhou et~al.(2018)Zhou, Rangarajan, and Gader]{zhou2018}
Yuan Zhou, Anand Rangarajan, and Paul~D. Gader.
\newblock A {G}aussian mixture model representation of endmember variability in
  hyperspectral unmixing.
\newblock \emph{IEEE Transactions on Image Processing}, 27\penalty0
  (5):\penalty0 2242--2256, 2018.

\end{thebibliography}
